\documentclass[11pt]{article}

\usepackage{color}
\usepackage{amsmath,amsfonts,amssymb,amsbsy}
\usepackage[english]{babel}
\usepackage{epstopdf}
\usepackage{bm}
\usepackage{mathtext}
\usepackage{amssymb}
\usepackage{graphics}
\usepackage[dvips]{graphicx}
\usepackage[left=2cm,right=2cm,top=2cm,bottom=2cm,bindingoffset=0cm]{geometry}
\usepackage{epstopdf}
\usepackage{soul}
\usepackage{dsfont}

\newtheorem{theorem}{\bf Theorem}[section]
\newtheorem{definition}{\bf Definition}[section]
\newtheorem{lemma}{\bf Lemma}[section]
\newtheorem{remark}{\bf Remark}[section]

\begin{document}

\begin{center}
{\bf \Large Regular approximate factorization of a class of matrix-function\\ with an unstable set of partial indices}
\end{center}

\begin{center}
G. Mishuris$^{1}$, S. Rogosin$^{2}$\\
$^{1}$Aberystwyth University, Penglais, SY23 3BZ Aberystwyth, UK\\
$^{2}$Belarusian State University, Nezavisimosti Ave., 4, 220030 Minsk, Belarus
\end{center}

\begin{abstract}
From the classic work of Gohberg and Krein (1958), it is well
known that the set of partial indices of a non-singular matrix
function may change depending on the properties of the original
matrix. More precisely, it was shown that if the difference
between the larger and the smaller partial indices is larger than
unity then, in any neighborhood of the original matrix function,
there exists another matrix function possessing a different set of
partial indices. As a result, the factorization of matrix
functions, being an extremely difficult process itself even in the
case of the canonical factorization, remains unresolvable or even
questionable in the case of a non-stable set of partial indices.
Such a situation, in turn, has became an unavoidable obstacle to
the application of the factorization technique. This paper sets
out to answer a less ambitious question than that of effective
factorizing matrix functions with non-stable sets of partial
indices, and instead focuses on determining the conditions which,
when imposed on the original matrix function, allow to construct
another matrix function that exhibits the same partial indices and
is close to the original matrix function.
\end{abstract}


\section{Introduction}

A given invertible square matrix $G \in\left( {\mathcal
C}({\mathbb R})\right)^{n\times n}$ is called factorizable if it
can be represented in the form
\begin{equation}
\label{fact}
G(x) = G^{-}(x) \Lambda(x)  G^{+}(x),
\end{equation}
with continuous invertible factors $G^{\pm}(x), (G^{\pm})^{-1}(x)$, possessing analytic continuation
into the corresponding half-plane $\Pi^{\pm}
= \{z = x + i y: {\mathrm{Im}}\, \pm z < 0\}$,  and



\begin{equation}
\label{fact1}
\Lambda(x) = {\mathrm{diag}}\, \left(\left(\frac{x - i}{x + i}\right)^{{\varkappa}_1}, \ldots, \left(\frac{x - i}{x + i}\right)^{{\varkappa}_n}\right),\; {\varkappa}_1, \ldots, {\varkappa}_n\in {\mathbb Z}.
\end{equation}
The representation (\ref{fact}) is called the {\it right (or
right-sided) factorization} on the real axis. It is also referred
to as right continuous factorization. If we have the
representation
\begin{equation}
\label{main_factorisation}
G(x) = G^{+}(x) \Lambda(x)  G^{-}(x),
\end{equation}
then it is called the {\it left (or left-sided) factorization}.
{If the right- (or left-)} factorization exists, then the
{integers} ${\varkappa}_1, \ldots, {\varkappa}_n$, called the {\it
partial indices}, are determined uniquely up to their order, e.g.
${\varkappa}_1 \geq \ldots \geq {\varkappa}_n$. The factors
$G^{-}$, $G^{+}$ are not unique (see, e.g., \cite{LitSpi87}). In
general, the partial indices for the right-factorization and for
the left-factorization are not necessarily the same. Throughout
{this paper we will only deal with the} right factorization. The
case of the left factorization can be handled analogously.

Several of the general facts on factorization have been presented
in the survey paper \cite{MishRog15a} (see also
\cite{GoKaSp03,LawAbr07}). In particular, it is well known that
the sum of the partial indices  is equal to the winding number (or
the Cauchy index) of the determinant of the given invertible
square matrix $G$:
\begin{equation}
\label{sum}
\sum_{j=1}^n{\varkappa}_j = {\varkappa}_{G} =
wind_{\mathbb R} \; det\, G(x) = \frac{1}{2\pi} \int\limits_{-\infty}^{+\infty} d (arg\, det\, G(x)).
\end{equation}

The factorization is called {the} {\it canonical factorization} if
all the partial indices are equal to 0, i.e. ${\varkappa}_1 =
\ldots = {\varkappa}_n = 0$.

It is said (see, e.g. \cite[p. 50]{GohKre58}) that a non-singular
matrix function $G(x)$ has a {\it stable set of partial indices}
if there exists $\delta > 0$ such that any matrix function $F(x)$
from the $\delta$-neighbourhood of $G(x)$ (i.e. $\|F - G\| <
\delta$) has the same set of partial indices (right or left). If
not, then $G(x)$ has {an} {\it unstable set of partial indices}.
It has been shown (see \cite{Boj58,GohKre58,GohKre58b} and also
\cite{Lit67}) that a set of partial indices ${\varkappa}_1,
\ldots, {\varkappa}_n$ is stable if and only if ${\varkappa}_1 -
{\varkappa}_n \leq 1$.

In the unstable case, a small deformation of the matrix function
$G(x)$ {\it can} lead to changes in the partial indices. More
precisely, there exists a sequence of matrix functions $F_k(x)$
that converge to $G(x)$, but which has a distinct set of partial
indices, i.e.
\begin{equation}
\label{sequence}
F_k(x) = F_k^{-}(x) \Lambda_A(x) F_k^{+}(x),
\end{equation}
{where} $\Lambda_A(x)\not= \Lambda(x)$.

We note that, in all the known examples illustrating such a
situation (see, e.g. \cite{Boj58,GohKre58}), the sequences of the
factors $F_k^{-}(x)$, $F_k^{+}(x)$ {\it do not} possess limiting
values (as $k\to\infty$) from the same chosen space as $G(x)$. On
the contrary, even in the case of unstable partial indices, we can
easily construct a sequence of factors $F_k^{-}(x)$, $F_k^{+}(x)$
in (\ref{sequence}) with the same partial indices as the original
matrix, i.e. $\Lambda_A(x)= \Lambda(x)$. Indeed, {in a simple
example we present the} pair $F_k^{-}(x) = G^{-}(x) +
\varepsilon_k H^{-}_k(x)$,  $F_k^{+}(x) = G^{+}(x) + \varepsilon_k
H^{+}_k(x)${,} where $H^{\mp}_k$ are arbitrary matrices belonging
to the same space as $G^{\mp}$, {such} that $\|H^\mp_k\|\le h_0$
and $\varepsilon_k\rightarrow 0$ as $k\rightarrow \infty$.

Let us now consider a matrix function $G\in {\mathcal G}{\mathcal
H}_{\mu}({\mathbb R})^{n\times n}$ that possesses a factorization.
Any matrix $G_{\varepsilon}\in {\mathcal G}{\mathcal
H}_{\mu}({\mathbb R})^{n\times n}$ that satisfies the following
asymptotic relation
$$
\|G_{\varepsilon}(x) - G(x)\| = O(\varepsilon),\;\; \varepsilon \rightarrow 0,
$$
will be called {a} {\it perturbation} of the matrix $G$.
\begin{definition}
\label{regular_per} Let $G\in {\mathcal G}{\mathcal
H}_{\mu}({\mathbb R})^{n\times n}$ be a given factorizable matrix.
Its perturbation $G_{\varepsilon}$ is considered {\it `regular'},
if there exists $\varepsilon_0 > 0$ such that the matrix
$G_{\varepsilon}$ possesses a bounded factorization (i.e.
$|G_{\varepsilon}^{\pm}(z)|\leq M$ for all $\varepsilon\in [0,
\varepsilon_0)$ and $z\in \Pi^{\pm}$). Otherwise the perturbation
is considered {\it `singular'}.
\end{definition}
\begin{lemma}
\label{pat_ind_per} The partial indices of the regular
perturbation $G_{\varepsilon}$ are the same as those of $G$.
\end{lemma}

{\bf Proof}
Let
$$
G(x) = G^{-}(x) \Lambda(x) G^{+}(x),
$$
and $G_{\varepsilon}$ be a regular perturbation of $G(x)$,
$$
G_{\varepsilon}(x) = G_{\varepsilon}^{-}(x) \Lambda_A(x) G_{\varepsilon}^{+}(x) =
G^{-}(x) \Lambda(x) G^{+}(x) +  O(\varepsilon).
$$
Hence
$$
\Lambda(x) = \left(G^{-}(x)\right)^{-1} G_{\varepsilon}^{-}(x) \Lambda_A(x) G_{\varepsilon}^{+}(x) \left(G^{+}(x)\right)^{-1} +  O(\varepsilon).
$$
By taking the limit as $\varepsilon\to +0$, we obtain $\Lambda_A =
\Lambda$ due to the uniqueness of the partial indices.

\begin{remark}
\label{stable_per} If the partial indices of $G$ satisfy the
condition ${\varkappa}_{max} - {\varkappa}_{min} \leq 1$, then any
perturbation of $G$ is regular. \end{remark}

\begin{remark}
\label{Ref_2}
To highlight the role of the condition of boundedness of the
factors{,} we present here a variant of the classical example {of}
Gohberg and Krein \cite[p. 264]{GohKre58}. Let us consider the
following matrix
\begin{equation}
\label{GK0}
G_{0}(x) = \left(\begin{array}{cc} \frac{x - i}{x + i} & 0 \\
0 & \frac{x + i}{x - i} \end{array}\right).
\end{equation}
It is clear that this matrix possesses {a} factorization{, where}
$G^{\pm}_{0}(x) = I$, $\Lambda(x) = G_{0}(x)$ and {the} partial
indices ${\varkappa}_1 = 1, {\varkappa}_2 = -1$. Consider a slight
perturbation of the matrix $G_{0}(x)$:
\begin{equation}
\label{GK1}
G_{\varepsilon}(x) = \left(\begin{array}{cc} \frac{x - i}{x + i} & \varepsilon \\
0 & \frac{x + i}{x - i} \end{array}\right),\;\;\; \varepsilon > 0.
\end{equation}
{We note} that for sufficiently small $\varepsilon > 0${, the}
matrices are close to each other, {such that}
$$
\|G_0(x) - G_{\varepsilon}(x)\| \leq \varepsilon.
$$
From the other hand, for all fixed $\varepsilon > 0$, the matrix
$G_{\varepsilon}(x)$ possesses the following factorization with
partial indices ${\varkappa}_1 = {\varkappa}_2 = 0$:
\begin{equation}
\label{GK2}
G_{\varepsilon}(x) = \left(\begin{array}{cc} 1 & 0 \\ \frac{1}{\varepsilon} \frac{x + i}{x - i} & 1
\end{array}\right) \cdot I \cdot \left(\begin{array}{cc} \frac{x - i}{x + i} & \varepsilon \\ - \frac{1}{\varepsilon}  & 0
\end{array}\right).
\end{equation}
For each fixed $\varepsilon > 0$, the factors in this
factorization admit analytic continuations into the corresponding
half-plane, where they are bounded. However, these factors are not
uniformly bounded for $\varepsilon\in [0, \varepsilon_0)$ for any
$\varepsilon_0 > 0$. Hence $G_{\varepsilon}(x)$ is a singular
perturbation of the above diagonal matrix $G_0(x)$ (we denote it
{by} $G_{\varepsilon}^{(s)}(x)$).
\end{remark}

\begin{figure}[h!]
\hspace*{-50mm}
\includegraphics [scale=0.85]{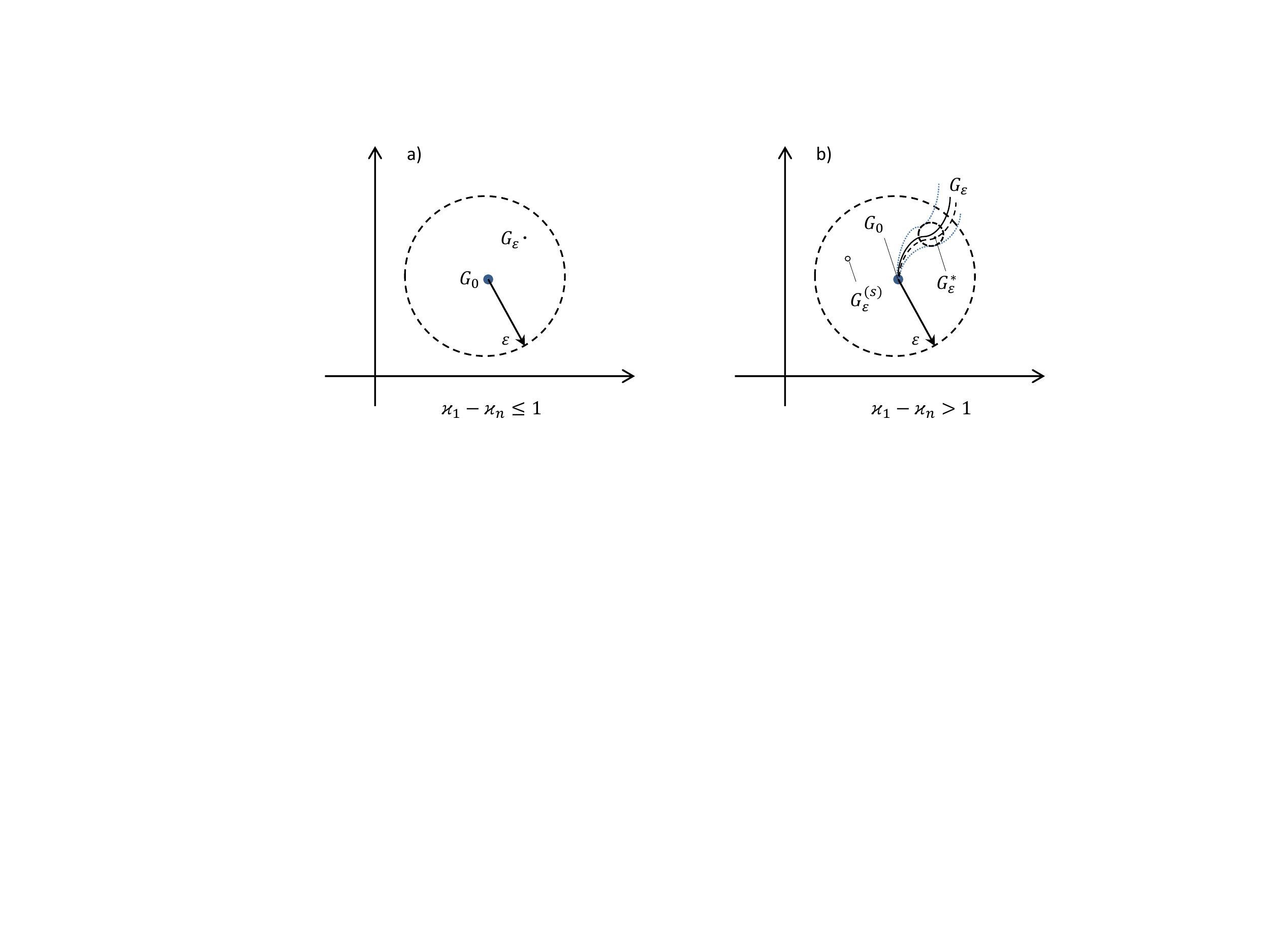}
\vspace*{-90mm}

    \caption{\footnotesize{Possible types of perturbations, $G_\varepsilon$,
    ($\|G_\varepsilon-G_0\|<\varepsilon$){, in the cases of} stable - a) and unstable - b)
    sets of partial indices of the matrix-function $G_0$.}}
\label{f0}
\end{figure}

This example provides a simple, but not unique, method  for the
construction of the singular perturbation
$G_{\varepsilon}^{(s)}(x)$ of any $n\times n$ diagonal matrix
$\Lambda(x) = {\mathrm{diag}}
\left\{\left(\frac{x-i}{x+i}\right)^{\varkappa_1},\ldots,\left(\frac{x-i}{x+i}\right)^{\varkappa_n}\right\}$.
Moreover, by replacing $\varepsilon$ with $\varepsilon^k$ in this
procedure, we can construct a singular perturbation for any
factorizable matrix $G(x) = G^{-}(x) \Lambda(x) G^{+}(x)$
(see Figure\, \ref{f0}),
\begin{equation}
\label{singular}
G_{\varepsilon}^{(s)}(x) = G^{-}(x) \Lambda_{\varepsilon}^{(s)}(x) G^{+}(x),
\end{equation}
which is arbitrarily close to the given matrix $G(x)$.
\begin{definition}
\label{guided} Let $G_{\varepsilon}(x)$ be a perturbation of a
factorizable matrix function $G_0(x)$ ($\varepsilon=0$). If there exists another
perturbation $G_{\varepsilon}^{\ast}(x)$ satisfying
$$
\|G_{\varepsilon}(x) - G_{\varepsilon}^{\ast}(x)\| = O(\varepsilon^k),\;\; \varepsilon \rightarrow 0,
$$
then we say that $G_{\varepsilon}^{\ast}$ is a ``k- {\it guided
perturbation''} of $G_{\varepsilon}$.
\end{definition}
It follows from (\ref{singular}) that for each regular perturbation $G_{\varepsilon}(x)$ of the
matrix function $G_0(x)$ there exists a singular k-guided perturbation $G_{\varepsilon}^{(s)}(x)$ for any $k\geq 1$.

The above mentioned properties are illustrated in Figure \ref{f0}.
In sub-figure a), we show that for any factorizable matrix
$G_0(x)$, with a stable set of partial indices, there is a
$\varepsilon$-neighbourhood containing only regular perturbations.
Sub-figure b) illustrates the case of unstable partial indices of
$G_0(x)$. Here, the situation is more delicate, as we can see that
in each $\varepsilon$-neighbourhood of $G_0(x)$ we can find either
regular or singular perturbations.

The aim of this paper is to consider the construction of a regular
k-guided ($k > 1$) perturbation $G_{\varepsilon}^{\ast}(x)$ for a
given perturbation $G_{\varepsilon}(x)$ of the matrix function
$G_0(x)$, with a known factorization {and unstable partial
indices}. For $k = 1$, this is trivial but gives no practical use.

The factorization technique is a powerful tool {in solving}
practical problems
\cite{AbrLaw95,AbrDavSmi08,Hur76,HurLun81,MishMovMov10,PiMiMo07,PiMiMo09,PiMiMo12}),
and any approximate factorization will allow a wide range of
practical problems to be tackled with some level of accuracy. The
establishment of an approximate (see e.g. \cite{Abr00}) or an
asymptotic procedure (see
\cite{Cri01,Kis13,Kis15,MishRog14a,MishRog15b}) is a challenging
 problem, since an exact factorization is possible only in a number
of special cases (see \cite{MishRog15a} and references therein).
Similarly, the mentioned non-uniqueness of the factorization
problem does not prevent it being effectively used in practice.
However one needs to be careful in using the approximate
factorization in the case of unstable indices, as it may introduce
not only quantitative but also qualitative deviation of the
approximate solution from the original one. Here, any links
between the partial indices of the factorization problem and the
particular physical properties of the problem in question are
crucial. This question is beyond the scope of this paper.

Below, we discuss whether, and under which conditions, it is
possible to find an $n\times n$ matrix function
$G_{\varepsilon}^{\ast}(x), x\in {\mathbb R}$, sufficiently close
to a given regular perturbation $G_{\varepsilon}(x)$ of the matrix
function $G_{0}(x)$, and possessing an unstable set of partial
indices. More exactly, we ask when it is possible to find
$G_{\varepsilon}^{\ast}(x)$ while preserving the partial indices
of $G_{\varepsilon}(x)$? To reach an answer to this question in
the case of unstable partial indices, a new definition of the
asymptotic factorization is given and applied. The method, as
proposed in \cite{MishRog14a}, \cite{MishRog15b}, is generalized
and employed. We find conditions under which our asymptotic
procedure is effective, and its properties and details are
illustrated by examples. The efficiency of the procedure is also
illuminated by numerical results.

The paper is organized as follows: In Sec.\, \ref{class} we
present necessary definitions and notations supplied by necessary
basic facts from factorization theory. Next, (Sec.\,
\ref{preserve}) we consider certain perturbations of matrices
factorized with unstable partial indices, and describe an
algorithm for the construction of an approximate factorization of
the perturbed matrices, while preserving the initial partial
indices. The conditions for realization of this algorithm are also
derived here, which are simply {\it solvability conditions} for a
certain boundary value problem. We also provide examples where the
solvability conditions are satisfied and unsatisfied. We conclude
with illustrations of the obtained numerical results and a
discussion of their importance in Sec.\, \ref{numerics}.

\section{Asymptotic factorization. Definitions}
\label{class}

To proceed, we make some necessary definitions. We denote by
${\mathcal H}_{\mu}({\mathbb R})^{n\times n}$, $n \in \mathbb{N}$,
the set of bounded matrix functions with locally
H\"older-continuous entries, endowed with the norm
$\|\cdot\|_{\mu}$:
\begin{equation}
\label{space} {\mathcal H}_{\mu}({\mathbb R})^{n\times n} =
\left\{G = (g_{ij}): {\mathbb R} \rightarrow {\mathcal M}^{n\times n}:  \|G\|_{\mu} = \max\limits_{1\leq i, j\leq n}
\|g_{ij}\|_{\mu} < \infty\right\}.
\end{equation}
In this article we consider only matrices of the class ${\mathcal
G} {\mathcal H}_{\mu}({\mathbb R})^{n\times n}$, where ${\mathcal
G}$ refers to {\it invertible} matrices. It should be noted that
our method can be also be applied to a wider class of matrix
functions.

Below, we give a definition of the asymptotic factorization only
in the case of the regular perturbation of a given matrix
function, since we lack for the moment a formal procedure which
distinguishes between the cases of regular and singular
perturbations and as yet have no useful example of the
construction of the asymptotic factorization of a singularly
perturbed matrix function.

\begin{definition}
\label{def} Let $G_0(x)\in {\mathcal G} {\mathcal
H}_{\mu}({\mathbb R})^{n\times n}$ be a given factorizable matrix
($G_0(x) = G_0^{-}(x) \Lambda(x) G_0^{+}(x)$) and
$G_{\varepsilon}(x)$ be its regular perturbation. We say that a
set of pairs of matrix functions, $ G^{-}_{\varepsilon,m}(x)$ and
$G^{+}_{\varepsilon,m}(x)$, ($m=1,2,\ldots N$), and a diagonal
matrix $\Lambda_A(x)$ of the form (\ref{fact1}) represent an
asymptotic factorization (of the order N) of the matrix function
$G(x)\in {\mathcal G} {\mathcal H}_{\mu}({\mathbb R})^{n\times n}$
 if the following conditions are satisfied:
\begin{enumerate}
\item
there exists a sequence of functions $\theta_k(\varepsilon)$,
$k=1,2,\ldots N+1${, that vanishes} at the point $\varepsilon=0${,
such that} for any $k=1,2,\ldots N$
\begin{equation}
\label{as_phi}
\theta_{k+1}(\varepsilon)=o\big(\theta_{k}(\varepsilon)\big),\quad
\varepsilon \to0;
\end{equation}
\item there exist matrices $ G^{\mp}_{\varepsilon,m}(x)$ of the form:
\begin{equation}
\label{as_G} G^{\mp}_{\varepsilon,m}(x)={G_{0}}^{\mp}(x)+\sum_{k=1}^m
\theta_k(\varepsilon)H^{\mp}_{\varepsilon,k}(x),
\end{equation}
\item
there exists $\varepsilon_0 > 0$ such that the matrices $H^{\mp}_{\varepsilon,k}(z)$
are  analytical in $\Pi^{\mp}$, respectively, and bounded in ${\mathcal H}_{\mu}({\mathbb
R})^{n\times n}$ uniformly with respect to $\varepsilon\in [0,\varepsilon_0)$,
\item
the following estimate is valid for any $m=1,2,\ldots N$
\begin{equation}
\label{as_full} G_0^{-}(x) \Lambda(x)
G_0^{+}(x)-G^{-}_{\varepsilon,m}(x) \Lambda_A(x)
G^{+}_{\varepsilon,m}(x)=O\big(\theta_{m+1}(\varepsilon)\big),\quad
\varepsilon \to0.
\end{equation}
\end{enumerate}
The representation
\begin{equation}
\label{asy_reg} G_{\varepsilon}^{\ast}(x) =
G^{-}_{\varepsilon,N}(x) \Lambda_A(x) G^{+}_{\varepsilon,N}(x),
\end{equation}
is called an asymptotic factorization  (of order $N$) of the
matrix $G_{\varepsilon}(x)$.
\end{definition}

We note that conditions \eqref{as_phi} and \eqref{as_G} guarantee
that the matrices ${G}^{\mp}_{\varepsilon,m}(z)$ and
$({G}^{\mp}_{\varepsilon,m}(z))^{-1}$ belong to the required
class, and thus both terms on the left-hand side of
\eqref{as_full} represent two essentially different
factorizations. As a simple example of the sequence
(\ref{as_phi}), we could consider
$\theta_k(\varepsilon)=\varepsilon^k$.

Some clarifications are required for this definition:
\begin{itemize}
\item
We are concerned with {\it regular} perturbations of the given
matrix, i.e. we are looking for {representations} (\ref{as_G})
possessing factors $G^{-}_{\varepsilon,m}(z)$ and
$G^{+}_{\varepsilon,m}(z)$ that are bounded in $z$ in the
corresponding half-planes, where our choice is motivated purely by
applications. In fact, we can replace the boundedness conditions
for $G^{-}_{\varepsilon,m}(z)$ and $G^{+}_{\varepsilon,m}(z)$, by
other conditions, such as polynomial growth/decay at infinity.

\item
The given definition does not demand uniqueness. Indeed, as was
demonstrated in the example in \cite{MishRog15b}, which
 consider the case of stable partial indices, there was no
uniqueness, even with the enforcement of additional conditions at
infinity.

\item
The parameter $N$ is also involved in the process of asymptotic
factorization. In the case when $N=\infty$ and the series is
converging, we can say that the asymptotic factorization becomes
the standard factorization, where the factors are defined by their
converging series.

\item
The method used  in \cite{MishRog15b}, in the case of stable
partial indices, allows to construct for some matrix functions the
factors of the factorization as converging asymptotic series, and
preserving the partial indices, i.e. $\Lambda_A(x)=\Lambda(x)$.
However, even in this case, no uniqueness can be guaranteed.

\item
The factors $G_{\varepsilon,m}^{\mp}(x)$ in the representation (\ref{as_G})
are continuous with respect to $\varepsilon\geq 0$.

\item
Although the asymptotic factorization is not unique, we can prove,
similarly to Lemma \ref{pat_ind_per}, the uniqueness of the
partial indices ($\Lambda_A(x)=\Lambda(x)$).

\item
If an asymptotic factorization of order $N > 1$ exists, then the
matrix function
$$
G_{\varepsilon}^{\ast}(x) = G^{-}_{\varepsilon,N}(z) \Lambda(x) G^{+}_{\varepsilon,N}(z)
$$
is an $(N+1)$-guided perturbation of the following perturbation
$G_{\varepsilon}(x)$ of the matrix $G(x)$:
$$
G_{\varepsilon}(x) = G^{-}_{\varepsilon,1}(z) \Lambda(x) G^{+}_{\varepsilon,1}(z).
$$
\end{itemize}

Although we only consider in this paper the factorization of a
matrix function on the real axis, we can tackle  in the same way
the factorization of matrices defined on any oriented curve
$\Gamma$ which divides the complex plane into two domains $D^{-}$
and $D^{+}$, by changing the diagonal entries in $\Lambda(x)$ to
$\left(\frac{x - t^{+}}{x - t^{-}}\right)^{{\varkappa}_j}$,
$t^{\mp}\in D^{\mp}$, or simply to $x^{{\varkappa}_j}$ if $0\in
D^{+}$.

Let us now consider a matrix function $G_{0}\in {\mathcal
G}{\mathcal H}_{\mu}({\mathbb R})^{n\times n}$, admitting a
factorization \eqref{fact}, with unstable partial indices
${\varkappa}_1, \ldots, {\varkappa}_n$ (${\varkappa}_1 \geq \ldots
\geq {\varkappa}_n$), and its perturbation $G_\varepsilon(x) \in
{\mathcal G}{\mathcal H}_{\mu}({\mathbb R})^{n\times n}$, which
depends on a small parameter $\varepsilon$ such that
\begin{equation}
\label{first_step_fac}
G_{\varepsilon}(x)\Bigl|_{\varepsilon=0} = G_{0}(x).
\end{equation}

Our motivating question is the following: how to distinguish a
class of possible perturbations that can be used to construct an
asymptotic procedure for the corresponding class of matrices,
according to the above definition \ref{def}, and how to find the
conditions under which this procedure can be realized? In the case
of stable partial indices, such a type of perturbation and the
corresponding asymptotic procedure was proposed in
\cite{MishRog14a,MishRog15b} for a class of matrix functions
related to certain applied problems. This procedure was shown to
be convergent for small values of the parameter $\varepsilon$.
Here, we extend the technique to the case of unstable partial
indices.

\section{Asymptotic factorization. Procedure}
\label{preserve}

Let us consider an invertible, bounded, locally H\"older continuous matrix
$G_\varepsilon(x) : {\mathbb R} \rightarrow {\mathcal G}{\mathcal H}_{\mu}({\mathbb R})^{n\times n}$ of the form
\begin{equation}
\label{Pr_new1} G_\varepsilon(x) = G_{0}(x) +
\theta_1(\varepsilon) N_{\varepsilon}(x),
\end{equation}
where $\theta_1(\varepsilon) ={o}(1)$ as $\varepsilon\to 0$ and
$N_{\varepsilon}$ is bounded and H\"older continuous on ${\mathbb
R}$. We suppose additionally that, when $\varepsilon = 0$, the
matrix $G_0(x)$ possesses a factorization with unstable partial
indices (${\varkappa}_1 - {\varkappa}_n \geq 2$), and has factors
$G_{0}^{\mp}(x)$ and $(G_{0}^{\mp}(x))^{\pm 1}$, which admit
analytic continuation into the semi-planes $\Pi^{\mp}$,
respectively, and which are bounded in $\overline{\Pi^{\mp}} =
\Pi^{\mp} \cup {\mathbb R}$.

We look for an asymptotic factorization of the matrix
{$G_\varepsilon(x)$} of the type (\ref{Pr_new1}, which is a
regular perturbation of $G_{0}(x)$, up to some stage of the
asymptotic procedure. For simplicity we will consider
$\theta_k(\varepsilon) = \varepsilon^k$ (see below Remark\,
\ref{choice_phi}, cf. \cite[Lemma 3.6]{MishRog14a}).

\subsection{First step of the asymptotic factorization}

First, we present the matrix {$G_\varepsilon(x)$} in the following
form:
\begin{equation}
\label{Pr_new6} G_\varepsilon(x) = \left(G_{0}^{-}(x) +
\varepsilon N_{1,\varepsilon}^{-}(x) (\Lambda^{+}(x))^{-1}\right)
\Lambda(x)\times
\end{equation}
$$
\left(G_{0}^{+}(x) + \varepsilon (\Lambda^{-}(x))^{-1}
N_{1,\varepsilon}^{+}(x)\right) + O(\varepsilon^2),
$$
where $\Lambda^{\mp}(x) = {\mathrm{diag}}\,
\left(\left(\frac{x-i}{x+i}\right)^{\varkappa_1^{\mp}}, \ldots{,}
\left(\frac{x-i}{x+i}\right)^{\varkappa_n^{\mp}}\right)$,
$\varkappa_j^{+} = \max \{\varkappa_j, 0\}, \varkappa_j^{-} =\max
\{0, -\varkappa_j\}$, and the unknown matrix-functions
$N_{1,\varepsilon}^{\mp}(x)$ must be be analytically extended into
$\Pi^{\mp}$, together with their inverses, and bounded in
$\overline{\Pi^{\mp}}$, respectively. Note that (\ref{Pr_new6})
differs from the representation used for the case of stable
partial indices (cf., \cite{MishRog14a}, \cite{MishRog15b}).

Comparing the term with parameter $\varepsilon$, we arrive at the
following boundary condition for $N_{1,\varepsilon}^{\mp}${:}
\begin{equation}
\label{Pr_new7} N_{1,\varepsilon}^{-}(x) \Lambda^{-}(x)
G_{0}^{+}(x) + G_{0}^{-}(x) \Lambda^{+}(x)
N_{1,\varepsilon}^{+}(x) = N_{\varepsilon}(x).
\end{equation}
For brevity, we introduce the following notation,
\begin{equation}
\label{Pr_new9} \widetilde{N_{1,\varepsilon}^{-}}(x) :=
\left(\tilde{n}_{1,ij}^{-}(z)\right)_{ij} = (G_{0}^{-}(x))^{-1}
N_{1,\varepsilon}^{-}(x),
\end{equation}
$$
\widetilde{N_{1,\varepsilon}^{+}}(x) :=
\left(\tilde{n}_{1,ij}^{+}(z)\right)_{ij} =
N_{1,\varepsilon}^{+}(x) (G_{0}^{+}(x))^{-1},
$$
$$
M_{0,\varepsilon}(x) := \left(m_{0,ij}(z)\right)_{ij} =
(G_{0}^{-}(x))^{-1} N_{\varepsilon}(x) (G_{0}^{+}(x))^{-1}.
$$
Hence, unknown matrix-functions
$\widetilde{N_{1,\varepsilon}^{\mp}}$ {have to} satisfy the
boundary condition
\begin{equation}
\label{Pr_new10} \widetilde{N_{1,\varepsilon}^{-}}(x)
\Lambda^{-}(x)  + \Lambda^{+}(x)
\widetilde{N_{1,\varepsilon}^{+}}(x) = M_{0,\varepsilon}(x),\;
x\in {\mathbb R}.
\end{equation}

In order to determine solvability conditions for this problem and
to find a representation for the solution, we present here a few
facts from the theory of boundary value problems (see
\cite{Gak77},\cite{Mus68}). It is known that any bounded locally
H\"older continuous function $f : {\mathbb R} \rightarrow {\mathbb
C}$ can be uniquely, up to arbitrary constant $c\in {\mathbb C}$,
represented as the sum of two functions which are analytic in
$\Pi^{-}$, $\Pi^{+}$, an bounded in $\overline{\Pi^{-}}$,
$\overline{\Pi^{+}}$, respectively,
\begin{equation}
\label{Pr_new11} f(x) = \left((\Omega_{0}^{-} f)(x) + c\right) +
\left((\Omega_{0}^{+} f)(x) - c\right),
\end{equation}
where $(\Omega_{0}^{\mp} f)(z)$ is the Cauchy type integral (see \cite[p. 52]{Gak77})
\begin{equation}
\label{Pr_new12} (\Omega_{0}^{\mp} f)(z) = \mp \frac{z - i}{2\pi
i} \int\limits_{-\infty}^{+\infty} \frac{f(\tau) d\tau}{(\tau -
i)(\tau - z)},\; z\in \Pi^{\mp}.
\end{equation}

Representation (\ref{Pr_new11}), and further formulas, remain
valid in the matrix case too.
Therefore, the {\it formal} solution to 
(\ref{Pr_new10}) has the following form
\begin{equation}
\label{Pr_new15} \widetilde{N_{1,\varepsilon}^{-}}(z) =
\left[(\Omega_{0}^{-} M_{0,\varepsilon})(z) + C^{0}\right]
(\Lambda^{-}(z))^{-1},
\end{equation}
\begin{equation}
\label{Pr_new16} \widetilde{N_{1,\varepsilon}^{+}}(z) =
(\Lambda^{+}(z))^{-1} \left[(\Omega_{0}^{+} M_{0,\varepsilon})(z)
- C^{0}\right],
\end{equation}
where $C^{0}$ is a constant $n\times n$ matrix, which is, in
coordinate-wise terms,
\begin{equation}
\label{Pr_new17} \left\{
\begin{array}{l}
\tilde{n}_{1,lj}^{-}(z) = (\Omega_{0}^{-} m_{0,lj})(z) + c_{lj}^{0},\; 1 \leq l \leq n, 1 \leq j \leq q, \\
\\
\tilde{n}_{1,lj}^{-}(z) = \left(\frac{z - i}{z + i}\right)^{-
{\varkappa}_{j}} \left[(\Omega_{0}^{-} m_{0,lj})(z) +
c_{lj}^{0}\right],\; 1 \leq l \leq n, q + 1 \leq j \leq n;
\end{array}
\right.
\end{equation}
\begin{equation}
\label{Pr_new18} \left\{
\begin{array}{l}
\tilde{n}_{1,lj}^{+}(z) = \left(\frac{z + i}{z - i}\right)^{{\varkappa}_{l}} \left[(\Omega_{0}^{+} m_{0,lj})(z) - c_{lj}^{0}\right],\; 1 \leq l \leq p, 1 \leq j \leq n, \\
\\
\tilde{n}_{1,lj}^{+}(z) = (\Omega_{0}^{+} m_{0,lj})(z) -
c_{lj}^{0},\; p + 1 \leq l \leq n, 1 \leq j \leq n.
\end{array}
\right.
\end{equation}
Such a representation gives the bounded solution 
of the first order asymptotic factorization problem
(\ref{Pr_new6}), with factors involving analytic matrices
$N_{1,\varepsilon}^{\mp}$ which are uniquely defined by
(\ref{Pr_new9}), if and only if certain solvability conditions are
satisfied. These conditions simply require that
$\widetilde{N_{1,\varepsilon}^{\mp}}(z)$ have no singular points
at $\mp i$. Partly, we can use arbitrary {constant} $c_{lj}^{0}$
(the entries of the matrix $C^{0}$), but not all the solvability
conditions are satisfied
by the proper choice of $c_{lj}^{0}$. 


\subsection{Solvability conditions}

Here, we present the necessary and sufficient solvability
conditions for boundary value problem (\ref{Pr_new10}), which is
equivalent to the first step of the asymptotic factorization (see
\cite[p. 120]{Gak77}).

\begin{itemize}
\item if for certain $k, q + 1 \leq k \leq n,$ we have ${\varkappa}_{k} = - 1$, then the boundedness of
$\widetilde{N_{1,\varepsilon}^{-}}(z)$ at $z = - i$ follows whenever we choose $c_{lk}^{0}$ such that
\begin{equation}
\label{sol_cond1} c_{lk}^{0} = \frac{1}{\pi}
\int\limits_{-\infty}^{+\infty} \frac{m_{0,lk}(\tau) d
\tau}{\tau^2 + 1},\; 1\leq l \leq n;
\end{equation}
\item if for certain $k, q + 1 \leq k \leq n,$ we have ${\varkappa}_{k} < - 1$, then the corresponding
$c_{lk}^{0}$ must be chosen as in (\ref{sol_cond1}), and the
entries $m_{0,lk}(\tau)$ have to satisfy conditions
\begin{equation}
\label{sol_cond2} \int\limits_{-\infty}^{+\infty}
\frac{m_{0,lk}(\tau) d \tau}{(\tau + i)^{r+1}} = 0,\; 1\leq r\leq
- {\varkappa}_{k} - 1,\; 1\leq l \leq n;
\end{equation}
\item if for certain $k, 1 \leq k \leq p,$ we have ${\varkappa}_{k} = 1$, then the boundedness of
$\widetilde{N_{1,\varepsilon}^{+}}(z)$ at
$z = i$ follows whenever we choose $c_{kj}^{0}$ such that
\begin{equation}
\label{sol_cond3} c_{kj}^{0} = 0,\; 1\leq j \leq n;
\end{equation}
\item if for certain $k, 1\leq k \leq p,$ we have ${\varkappa}_{k} > 1$, then the corresponding
$c_{kj}^{0}$ must be chosen as in (\ref{sol_cond3}), and the
entries $m_{0,kj}(\tau)$ have to satisfy conditions
\begin{equation}
\label{sol_cond4} \int\limits_{-\infty}^{+\infty}
\frac{m_{0,kj}(\tau) d \tau}{(\tau - i)^{r+1}} = 0,\; 1\leq r \leq
{\varkappa}_{k} - 1,\; 1\leq j \leq n;
\end{equation}
\item if the pair $(l,j)$ is such that $1\leq l \leq p, q + 1 \leq j \leq n$, then additional solvability
conditions must satisfy
\begin{equation}
\label{sol_cond5} \int\limits_{-\infty}^{+\infty}
\frac{m_{0,lj}(\tau) d \tau}{\tau^2 + 1} = 0,\; 1\leq l \leq p, q
+ 1 \leq j \leq n;
\end{equation}
\item if the pair $(l,j)$ is such that either $1\leq l\leq n, 1\leq j \leq q$, or
$p + 1\leq l\leq n, 1\leq j \leq n$, then we have no condition on
the entries $m_{0,lj}(\tau)$; the {corresponding constants
$c_{lj}^{0}$} can take arbitrary value.
\end{itemize}

\begin{theorem}
\label{first_order} Formula (\ref{Pr_new6}) gives the first order
bounded asymptotic factorization for all $\varepsilon$ smaller
than a certain positive $\varepsilon_1$ if ant only if the
solvability conditions (\ref{sol_cond2}), (\ref{sol_cond4}),
(\ref{sol_cond5}) are satisfied, and the constants $c_{ij}^{0}$
are chosen accordingly.
\end{theorem}

{\bf Proof}
If the conditions of the theorem are satisfied, then the matrix
functions $\widetilde{N_{1,\varepsilon}^{\mp}}(z)$ give a bounded
solution to the problem (\ref{Pr_new10}). Moreover, the matrices
$\widetilde{N_{1,\varepsilon}^{-}}(z) \Lambda^{-}(z)$,
$\Lambda^{+}(z) \widetilde{N_{1,\varepsilon}^{+}}(z)$ are bounded
in {the} neighbourhoods of $z = -i$, $z = i$, respectively. By
choosing sufficiently small $\varepsilon_1 > 0$, we can guarantee
that the matrix functions $G_{0}^{-}(z) + \varepsilon
N_{1,\varepsilon}^{-}(z) (\Lambda^{+}(z))^{-1}$, $G_{0}^{+}(z) +
\varepsilon (\Lambda^{-}(z))^{-1} N_{1,\varepsilon}^{+}(z)$ are
invertible in the corresponding semi-planes. Thus, for
$\varepsilon\in [0, \varepsilon_1)$, formula (\ref{Pr_new6}) gives
the first order bounded asymptotic factorization.

To demonstrate the necessity of the theorem's conditions, we
suppose that formula  (\ref{Pr_new6}) gives the first order
bounded asymptotic factorization. Then the matrix functions
${N_{1,\varepsilon}^{\mp}}(z)$ have to satisfy boundary condition
(\ref{Pr_new7}), being analytically extended into $\Pi^{\mp}$
together with their inverses, and bounded in
$\overline{\Pi^{\mp}}$, respectively. The boundary value problem
(\ref{Pr_new7}) is equivalent to (\ref{Pr_new10}), and for
invertibility of matrices  ${N_{1,\varepsilon}^{\mp}}(z)$ we must,
in particular, have boundedness of the matrix functions
(\ref{Pr_new15}), (\ref{Pr_new16}) in the neighbourhoods of $z = -
i$, $z = i$, respectively. The latter leads to the necessity of
the conditions of the theorem.

\begin{remark}
\label{sol_cond_number} The numbers of solvability conditions and
conditions on the choice of the constants satisfy the following
relations.
\begin{itemize}
\item The number of solvability conditions is given by
\begin{equation}
\label{sol_condl} \sum\limits_{j=q + 1}^{n} (- {\varkappa}_{j} - 1)n +
\sum\limits_{i=1}^{p} ({\varkappa}_{i} - 1)n + (n - q)p.
\end{equation}
\item $(n - q)n$ constants $c_{ij}^{0}$ are chosen according to (\ref{sol_cond1}) and $n p$ constants
$c_{ij}^{0}$ are equal to 0, as in (\ref{sol_cond3}). In the $(n -
q) p$ cases described in (\ref{sol_cond5}) these choices of the
constants  $c_{ij}^{0}$ must coincide.
\item $n(n - p + q)$ constants  $c_{ij}^{0}$ can be chosen arbitrarily.
\end{itemize}
\end{remark}
\begin{remark}
\label{manifold1} The obtained result can be interpreted in the
following manner. Let the matrix function $G_\varepsilon(x)$ be a
perturbation of $G_0(x)$. Then, in particular, $G_\varepsilon(x)$
is in the $\varepsilon$-neighbourhood of $G_{0}(x)$ (see Figure\,
\ref{f0} a)). If the matrix $G_\varepsilon(x)$ satisfies the above
solvability conditions, then there exists for all sufficiently
small $\varepsilon$ the matrix
\begin{equation}
\label{neigh1} G_{\varepsilon}^{\ast}(x) = \left(G_{0}^{-}(x) +
\varepsilon N_{1,\varepsilon}^{-}(x) (\Lambda^{+}(x))^{-1}\right)
\Lambda(x) \left(G_{0}^{+}(x) + \varepsilon (\Lambda^{-}(x))^{-1}
N_{1,\varepsilon}^{+}(x)\right),
\end{equation}
which possesses a factorization with the same unstable set of
partial indices as $G_0(x)$. The matrix
$G_{\varepsilon}^{\ast}(x)$ {is in the}
$\varepsilon^2$-neighbourhood of $G_\varepsilon(x)$ (see Figure
\ref{f0} b)). This means that for each point of linear manifold of
{the} matrices $G_\varepsilon(x)$, as defined by (\ref{Pr_new6}),
which satisfies the solvability conditions, there exists a point
(matrix $G_{\varepsilon}^{\ast}(x)$) in its
$\varepsilon^2$-neighbourhood which preserves the initial partial
indices, i.e. according to Definition\, \ref{guided}, the latter
matrix is the regular 2-guided perturbation.
\end{remark}

\subsection{Further steps of the asymptotic factorization}

Let the solvability conditions be satisfied and the constants
$c_{ij}^{0}$ chosen accordingly. By solving the corresponding
boundary value problems we can refine the first order
factorization up to the $r$-th step
of the factorization using the representation 
\begin{equation}
\label{PrF_new5} G_\varepsilon(x) = \left(G_{0}^{-}(x) +
\varepsilon N_{1,\varepsilon}^{-}(x) (\Lambda^{+}(x))^{-1} +
\ldots + \varepsilon^r N_{r,\varepsilon}^{-}(x)
(\Lambda^{+}(x))^{-1}\right) \Lambda(x) \times
\end{equation}
$$
\times \left(G_{0}^{+}(x) + \varepsilon (\Lambda^{-}(x))^{-1}
N_{1,\varepsilon}^{+}(x) + \ldots + \varepsilon^r
(\Lambda^{-}(x))^{-1} N_{r,\varepsilon}^{+}(x)\right) +
O(\varepsilon^{r+1}),
$$
which leads to the boundary value problem
\begin{equation}
\label{PrF_new6} \widetilde{N_{r,\varepsilon}^{-}}(x)
\Lambda^{-}(x)  + \Lambda^{+}(x)
\widetilde{N_{r,\varepsilon}^{+}}(x) = M_{r-1, \varepsilon }(x),\;
x\in {\mathbb R},
\end{equation}
\begin{equation}
\label{PrF_new7} \widetilde{N_{r,\varepsilon}^{-}}(x) :=
(G_{0}^{-}(x))^{-1} N_{r,\varepsilon}^{-}(x),\;
\widetilde{N_{r,\varepsilon}^{+}}(x) := N_{r,\varepsilon}^{+}(x)
(G_{0}^{+}(x))^{-1},
\end{equation}
$$
M_{r-1,\varepsilon}(x) := - (G_{0}^{-}(x))^{-1}
\left[N_{1,\varepsilon}^{-} (x) N_{r-1,\varepsilon}^{+}(x) +
\ldots + N_{r-1,\varepsilon}^{-} (x)
N_{1,\varepsilon}^{+}(x)\right] (G_{0}^{+}(x))^{-1}.
$$
The {\it formal} solution to problem (\ref{PrF_new6}) can be
{presented} as
\begin{equation}
\label{PrF_new8} \widetilde{N_{r,\varepsilon}^{-}}(z) =
\left[(\Omega_{0}^{-} M_{r-1,\varepsilon})(z) + C^{r-1}\right]
(\Lambda^{-}(z))^{-1},
\end{equation}
\begin{equation}
\label{PrF_new9} \widetilde{N_{r,\varepsilon}^{+}}(z) =
(\Lambda^{+}(z))^{-1} \left[(\Omega_{0}^{+}
M_{r-1,\varepsilon})(z) - C^{r-1}\right],
\end{equation}
which features a new constant matrix $C^{r-1}$. It becomes the
solution for the considered class if and only if the solvability
conditions (\ref{sol_cond1})-(\ref{sol_cond5}) are satisfied (in
this case we replace the functions $m_{0,ij}(x)$ with the
functions $m_{r-1,ij}(x)$, and the constants $c^{0}_{ij}$ with the
constants $c^{r-1}_{ij}$), while the constants $c^{r-1}_{ij}$ are
chosen accordingly.

If at a certain step $r = N+1$, at least one solvability condition
fails, then the procedure for the the asymptotic factorization is
stopped at this point.

\begin{remark}
\label{manifoldN} In this case, we summarize the situation as
follows. Let the matrix function $G_\varepsilon(x)$ be a regular
perturbation of $G_0(x)$. In particular, (as for $N=1$)
$G_\varepsilon(x)$ {is in} the $\varepsilon$-neighbourhood of
$G_0(x)$. If the matrix function $G_\varepsilon(x)$ satisfies the
above solvability conditions at each step $r, 1 \leq r \leq N,$
then for all sufficiently small $\varepsilon$ there exists a
matrix
\begin{equation}
\label{neigh1} G_{N,\varepsilon}^{\ast}(x) = (G_{0}^{-}(x) +
\sum\limits_{r=1}^{N} \varepsilon^r N_{r,\varepsilon}^{-}(x)
(\Lambda^{+}(x))^{-1}) \Lambda(x) (G_{0}^{+}(x) +
\sum\limits_{r=1}^{N} \varepsilon^r (\Lambda^{-}(x))^{-1}
N_{r,\varepsilon}^{+}(x)),
\end{equation}
which possesses a factorization with the same set of unstable
partial indices as $G_0(x)$. The matrix
$G_{N,\varepsilon}^{\ast}(x)$ is in the
$\varepsilon^{N+1}$-neighbourhood of $G_\varepsilon(x)$. This
means that for each point of linear manifold of the matrices
$G_\varepsilon(x)$ that satisfies the solvability conditions,
there exists an $(N+1)$-guided perturbation. Thus, with a larger
number of steps we can proceed in our approximate factorization
more closely to  the index-preserving approximation to a given
matrix $G_\varepsilon(x)$.
\end{remark}
\begin{remark} \label{failed} If at least one solvability
condition fails at some $N$-th step of the approximation, then we
may only construct an approximate factorization up to the order
$N-1$. If a solvability condition fails at the first step of
approximation, then we do not have a tool to construct a regular
$k$-guided perturbation for any $k >1$.
\end{remark}

\subsection{Example of the perturbed matrix satisfying the first order solvability conditions}
\label{Ex_preserve}


We apply the above described asymptotic procedure to the matrix
function $G_{\varepsilon}(x)$ {of the} form
\begin{equation}
\label{UnEx2_1} G_{\varepsilon}(x) = \left(
\begin{array}{cc}
 \frac{x^2 + x i (-18 + 8 e^{i  \varepsilon  x} + 8 e^{- i  \varepsilon  x}) - 1}{x^2 + 1} &
 \frac{x i (24 - 12 e^{i  \varepsilon  x} - 12 e^{- i  \varepsilon  x})}{x^2 + 1} \\
 & \\
\frac{x i (- 12 + 4 e^{i  \varepsilon  x} + 8 e^{- i  \varepsilon
x})}{x^2 + 1} & \frac{x^2 + x i (18 - 8 e^{i  \varepsilon  x} - 8
e^{- i  \varepsilon  x}) - 1}{x^2 + 1}
\end{array}
\right)
\end{equation}
and show that this matrix possesses an asymptotic factorization
with the same partial indices as $G_{0}(x)$. Here the matrix
function $G_{0}(x)$ is given by (\ref{GK0}).

The matrix function $G_{\varepsilon}(x)$ can be represented in the following form
\begin{equation}
\label{UnEx2_4} G_{\varepsilon}(x) = \Lambda(x) +
N_{\varepsilon}(x),
\end{equation}
where $\Lambda(x)$ is the same as in the Remark \ref{Ref_2} and
the matrix function $N_{\varepsilon}(x)$ is given by
\begin{equation}
\label{UnEx2_5} N_{\varepsilon}(x) = \left(
\begin{array}{cc}
 \frac{x i (-16 + 8 e^{i  \varepsilon  x} + 8 e^{- i  \varepsilon  x})}{x^2 + 1} &
 \frac{x i (24 - 12 e^{i  \varepsilon  x} - 12 e^{- i  \varepsilon  x})}{x^2 + 1} \\
 & \\
\frac{x i (- 12 + 4 e^{i  \varepsilon  x} + 8 e^{- i  \varepsilon
x})}{x^2 + 1} & \frac{x i (16 - 8 e^{i  \varepsilon  x} - 8 e^{- i
\varepsilon  x})}{x^2 + 1}
\end{array}
\right),
\end{equation}
or
$$
N_{\varepsilon}(x) = \left(
\begin{array}{cc}
 -\frac{32 x i \sin^2 \frac{ \varepsilon  x}{2}}{x^2 + 1} &
 \frac{48 x i \sin^2 \frac{ \varepsilon  x}{2}}{x^2 + 1} \\
 & \\
\frac{- 24 x i \sin \frac{ \varepsilon  x}{2} (\sin \frac{
\varepsilon  x}{2} - i \cos \frac{ \varepsilon  x}{2})}{x^2 + 1} &
\frac{32 x i \sin^2 \frac{ \varepsilon x}{2}}{x^2 + 1}
\end{array}
\right) = \frac{x \sin \frac{\varepsilon x}{2}}{x^2 + 1} \left(
\begin{array}{cc}
 -32 i \sin \frac{ \varepsilon  x}{2} &
 48 i \sin \frac{ \varepsilon  x}{2} \\
 & \\
- 24 i e^{-i \frac{\varepsilon  x}{2}} &
32 i \sin \frac{\varepsilon x}{2}
\end{array}
\right).
$$
Thus, $G_{\varepsilon}(x)$ can be thought of as a small
perturbation of the matrix function $G_{0}(x) = \Lambda(x)$
 ($G_0^\pm(x)=I$). The matrix
function $N_{\varepsilon}(x)$ takes the following representation
(uniform in $x\in {\mathbb R}$ and in $ \varepsilon $) on any
finite interval:
$$
N_{\varepsilon}(x) = \phi(x,\varepsilon) {\widetilde
N}_{\varepsilon}(x),\quad \phi(x,\varepsilon) =  \frac{x \sin \frac{\varepsilon x}{2}}{x^2 + 1},
$$
where ${\widetilde N}_{\varepsilon}(x)$ is a bounded matrix.
\begin{remark}
\label{choice_phi} The introduced small parameter
$\phi(x,\varepsilon)$ has the following properties (cf.
\cite[Lemma 3.6]{MishRog14a}):
\begin{equation}
\label{choice_phi1}
\phi(x,\varepsilon) = O(\varepsilon),\;\;\; \forall 0 < \varepsilon < \varepsilon_0,
\end{equation}
\begin{equation}
\label{choice_phi2}
\phi(x,\varepsilon) = O\left(\frac{1}{x}\right),\;\;\; |x|\to +\infty.
\end{equation}
We note here that $\theta_1(\varepsilon) =
\max\limits_{x\in\overline{\mathbb R}} \phi(x,\varepsilon)$. In
our case, we can prove that $\theta_1(\varepsilon) =
O(\varepsilon)$. We can thus later use an artificial small
parameter $\varepsilon$ instead of $\theta_1(\varepsilon)$.
\end{remark}

\begin{remark}
\label{choice_phi1} In fact, the first order decay of
$\phi(x,\varepsilon)$ at infinity (\ref{choice_phi2}) is crucial
to the behaviour of $\theta_1(\varepsilon)$ with respect to
$\varepsilon$. If, for example, $\phi(x,\varepsilon) =
O\left(\frac{1}{\sqrt{x}}\right)$, then it leads to only
$\theta_1(\varepsilon) = O(\varepsilon^{1/2})$.
\end{remark}

$1^{0}.$ The first step of the asymptotic factorization procedure.

We look for a pair of matrix functions $N_{1, \varepsilon
}^{\pm}(x)$, which are an approximate solution, up to
$\varepsilon^1$, of the functional equation
\begin{equation}
\label{UnEx2_7} G_{ \varepsilon }(x) = \left(I + N_{1, \varepsilon }^{-}(x)
\left(\Lambda^{+}(x)\right)^{(-1)}\right) \Lambda(x) \left(I +
 \left(\Lambda^{-}(x)\right)^{(-1)} N_{1, \varepsilon }^{+}(x) \right) + O(\varepsilon^2).
\end{equation}
We remind here that $G_0^\pm(x)=I$ (cf. \ref{Pr_new7}). The
approximate solution to (\ref{UnEx2_7}) can be found from the
matrix boundary value problem (\ref{Pr_new10}) that takes in this
case the form:
\begin{equation}
\label{UnEx2_8} \Lambda^{+}(x) N_{1, \varepsilon }^{+}(x) + N_{1,
\varepsilon }^{-}(x) \Lambda^{-}(x) = N_{ \varepsilon }(x),
\end{equation}
where $ M_{0, \varepsilon }(x)= N_{ \varepsilon }(x)$.

Bounded solutions to (\ref{UnEx2_8}) have to satisfy the relation
\begin{equation}
\label{UnEx2_11} \Lambda^{+}(x) N_{1, \varepsilon }^{+}(x) = M_{0,
\varepsilon }^{+}(x) - C_0,\quad \Lambda^{-}(x) N_{1, \varepsilon
}^{-}(x) = M_{0, \varepsilon }^{-}(x) + C_0,
\end{equation}
where $C_0 = (c_{ij}^{0})$ is a constant matrix. Hence,
\begin{equation}
\label{UnEx2_12} N_{1, \varepsilon }^{+}(x) = \left(
\begin{array}{cc}
\frac{x + i}{x - i}(m_{0,11}^{+} - c_{11}^{0}) &
\frac{x + i}{x - i}(m_{0,12}^{+} - c_{12}^{0}) \\
\\
m_{0,21}^{+} - c_{21}^{0} & m_{0,22}^{+} - c_{22}^{0}
\end{array}
\right),
\end{equation}
\begin{equation}
\label{UnEx2_13} N_{1, \varepsilon }^{-}(x) = \left(
\begin{array}{cc}
m_{0,11}^{-} + c_{11}^{0} &
\frac{x - i}{x + i}(m_{0,12}^{-} + c_{12}^{0}) \\
\\
m_{0,21}^{-} + c_{21}^{0} & \frac{x - i}{x + i}(m_{0,22}^{-} + c_{22}^{0})
\end{array}
\right).
\end{equation}
For analyticity of $N_{1, \varepsilon }^{+}, N_{1,
\varepsilon}^{-}$ in the corresponding half-planes, it is
necessary and sufficient that the following conditions be
fulfilled,
\begin{itemize}
\item{} $c_{11}^{0} = m_{0,11}^{+}(i), \quad c_{22}^{0} = - m_{0,22}^{-}(-i)$;\\[2mm]
\item{} the constant $c_{21}^{0}$ is chosen {arbitrarily};\\[2mm]
\item{} {the} solvability condition $m_{0,12}^{+}(i) = -
m_{0,12}^{-}(-i)$ holds{;}\\[2mm]
\item{} $c_{12}^{0} =
m_{0,12}^{+}(i)$.
\end{itemize}

In the case of the matrix function $N_{ \varepsilon }(x)$ given by
(\ref{UnEx2_5}), we have
\begin{equation}
\label{UnEx2_14} N_{1, \varepsilon }^{+}(x) = \left(
\begin{array}{cc}
 \frac{x + i}{x - i}\left(\frac{-8i(1 - e^{- \varepsilon })}{x + i} + \frac{8 x i (e^{i  \varepsilon  x} - e^{-  \varepsilon })}{x^2 + 1} - c_{11}^{0}\right) &
\frac{x + i}{x - i}\left(\frac{12i(1 - e^{- \varepsilon })}{x + i} - \frac{12 x i (e^{i  \varepsilon  x} - e^{-  \varepsilon })}{x^2 + 1} - c_{12}^{0}\right) \\
 & \\
\frac{-6i(1 - e^{- \varepsilon })}{x + i} + \frac{4 x i
(e^{i  \varepsilon  x} - e^{-  \varepsilon })}{x^2 + 1} -
c_{21}^{0} & \frac{8i(1 - e^{- \varepsilon })}{x + i} -
\frac{8 x i (e^{i  \varepsilon  x} - e^{-  \varepsilon })}{x^2 +
1} - c_{22}^{0}
\end{array}
\right),
\end{equation}
\begin{equation}
\label{UnEx2_15} N_{1, \varepsilon }^{-}(x) = \left(
\begin{array}{cc}
 \frac{-8i(1 - e^{- \varepsilon })}{x - i} + \frac{8 x i (e^{-i  \varepsilon  x} - e^{-  \varepsilon })}{x^2 + 1} + c_{11}^{0} &
\frac{x - i}{x + i}\left(\frac{12i(1 - e^{- \varepsilon })}{x - i} - \frac{12 x i (e^{-i  \varepsilon  x} - e^{-  \varepsilon })}{x^2 + 1} + c_{12}^{0}\right) \\
 & \\
\frac{-6i(1 - e^{- \varepsilon })}{x - i} + \frac{8 x i
(e^{-i  \varepsilon x} - e^{-  \varepsilon })}{x^2 + 1} +
c_{21}^{0} & \frac{x - i}{x + i}\left(\frac{8i(1 - e^{-
\varepsilon })}{x - i} - \frac{8 x i (e^{-i  \varepsilon x} - e^{-
\varepsilon })}{x^2 + 1} + c_{22}^{0}\right)
\end{array}
\right).
\end{equation}
Here
\begin{equation}
\label{UnEx2_16} c_{11}^{0} = m_{0,11}^{+}(i) = - 4 (1 - e^{-
\varepsilon }) - 4 \varepsilon  e^{- \varepsilon }, \quad c_{22}^{0} =
- m_{0,22}^{-}(-i) = 4 (1 - e^{- \varepsilon }) + 4 \varepsilon
e^{- \varepsilon },
\end{equation}
the solvability condition is satisfied
\begin{equation}
\label{UnEx2_17} m_{0,12}^{+}(i) = - m_{0,12}^{-}(-i)  = 6 (1 -
e^{- \varepsilon }) + 6 \varepsilon  e^{- \varepsilon },
\end{equation}
and thus the constant $c_{12}^{0}$ can be chosen accordingly
\begin{equation}
\label{UnEx2_18} c_{12}^{0} = m_{0,12}^{+}(i) = 6 (1 - e^{-
\varepsilon }) + 6
 \varepsilon  e^{- \varepsilon }.
\end{equation}
Finally, the constant $c_{21}^{0}$ can be chosen arbitrarily.

Thus, the first order approximation $G_{1,\varepsilon}^{\ast}(x)$
for the factorization of $G_{\varepsilon}(x)$  is given by the
following formula
\begin{equation}
\label{first_app}
G_{1,\varepsilon}^{\ast}(x) := \left(I + N_{1, \varepsilon }^{-}(x)
\left(\Lambda^{+}(x)\right)^{(-1)}\right) \Lambda(x) \left(I +
 \left(\Lambda^{-}(x)\right)^{(-1)} N_{1, \varepsilon }^{+}(x) \right),
\end{equation}
where matrices $N_{1, \varepsilon }^{\pm}(x)$ are presented in
(\ref{UnEx2_14}), (\ref{UnEx2_15}) with the above described choice
of constants.

In order to estimate the quality of the approximation, it is
customary to define the following remainder matrix
\begin{equation}
\label{difference}
\Delta K_{1,\varepsilon}(x) := G_{\varepsilon}(x) - G_{1,\varepsilon}^{\ast}(x).
\end{equation}
Direct calculations show that $\Delta K_{1,\varepsilon}(x) =
O(\varepsilon^2)$ as $\varepsilon\to +0$ and thus
$G_{1,\varepsilon}^{\ast}(x)$ is the 2-guided perturbation for the
matrix $G_{\varepsilon}(x)$.

\begin{remark}
\label{infty} Matrix $\Delta K_{1,\varepsilon}(x)$ has an
interesting behaviour, as a consequence of a special property of
the matrix $N_{\varepsilon}(x): \; n_{11} = - n_{22}$. Namely, it
tends to
 the diagonal matrix as $x\to \infty$, specifically,
\begin{equation}
\label{difference_infty}
\Delta K_{1,\varepsilon}(\infty) =  \left(\begin{array}{cc} (c_{11}^{0})^2 + c_{12}^{0} c_{21}^{0} & 0 \\
0 & (c_{22}^{0})^2 + c_{12}^{0} c_{21}^{0} \end{array}\right).
\end{equation}
Thus, by taking $c_{21}^{0} = 0$ we have
$$
\Delta K_{1,\varepsilon}(\infty) =  16 (1 - e^{-\varepsilon} + \varepsilon e^{-\varepsilon})^2 I,
$$
and by taking $c_{21}^{0} =  - 8/3 (1 - e^{-\varepsilon} + \varepsilon e^{-\varepsilon})$
we have
$$
\Delta K_{1,\varepsilon}(\infty) = \left(\begin{array}{cc} 0 & 0 \\
0 & 0 \end{array}\right).
$$
\end{remark}
The above two characteristic values of the constant $c_{21}^{0}$
will be used in our numerical description of the behaviour of the
remainder $\Delta K_{1,\varepsilon}(x)$ of the first order
approximate factorization of the matrix (\ref{UnEx2_1}).

In this example we have restricted our calculation to only the
first step of the approximation. In principle, the procedure for
the next steps has been already been described. However, there is
no guarantee that the next step will be successful and a higher
order guided perturbation will have been derived.

\subsection{Example of a matrix which does not satisfy the solvability conditions}

Simple changes to the matrix ${G}_{\varepsilon}(x)$ can lead to a
violation of the solvability conditions for the corresponding
boundary value problem. Let us consider
\begin{equation}
\label{UnEx3_1} \hat{G}_{ \varepsilon }(x) = \left(
\begin{array}{cc}
 \frac{x^2 + x i (-18 + 8 e^{i  \varepsilon  x} + 8 e^{- i  \varepsilon  x}) - 1}{x^2 + 1} &
 \frac{x i (24 - 16 e^{i  \varepsilon  x} - 8 e^{- i  \varepsilon  x})}{x^2 + 1} \\
 & \\
\frac{x i (- 12 + 4 e^{i  \varepsilon  x} + 8 e^{- i  \varepsilon
x})}{x^2 + 1} & \frac{x^2 + x i (18 - 8 e^{i  \varepsilon  x} - 8
e^{- i  \varepsilon  x}) - 1}{x^2 + 1}
\end{array}
\right).
\end{equation}
As before, $\hat{G}_{0}(x) = {G}_{0}(x)$, and thus
$\hat{G}_{0}(x)$ possesses a factorization with partial indices
${\varkappa}_1 = 1, {\varkappa}_2 = -1$.

We note that
$$
\hat{G}_{\varepsilon }(x) = {G}_{ \varepsilon }(x) +  \left(
\begin{array}{cc}
 0 &
 -\frac{x i (4 e^{i  \varepsilon  x} - 4 e^{- i  \varepsilon  x})}{x^2 + 1} \\
0 & 0
\end{array}
\right).
$$

We apply the above described asymptotic procedure to our matrix
$\hat{G}_{\varepsilon}(x)$, and show that this matrix cannot
possess a bounded first order asymptotic factorization with the
same partial indices as $\hat{G}_{0}(x)$.

The corresponding matrix is given by
\begin{equation}
\label{UnEx3_4} \hat{N}_{\varepsilon}(x) :=
\hat{G}_{\varepsilon}(x) - \Lambda(x) = \left(
\begin{array}{cc}
\frac{x i (- 16 + 8 e^{i  \varepsilon  x} + 8 e^{- i  \varepsilon  x})}{x^2 + 1} & \frac{x i (24 - 16 e^{i  \varepsilon  x} - 8 e^{- i  \varepsilon  x})}{x^2 + 1}\\
\frac{x i (- 12 + 4 e^{i  \varepsilon  x} + 8 e^{- i  \varepsilon
x})}{x^2 + 1} & \frac{x i (16 - 8 e^{i  \varepsilon  x} - 8 e^{- i
\varepsilon  x})}{x^2 + 1}
\end{array}
\right) =
\end{equation}
$$
= \frac{x \sin \frac{\varepsilon x}{2}}{x^2 + 1} \left(
\begin{array}{cc}
 -32 i \sin \frac{ \varepsilon  x}{2} &
 2 i (24 \sin \frac{ \varepsilon  x}{2} + i \cos \frac{ \varepsilon  x}{2})\\
 & \\
- 24 i (\sin \frac{ \varepsilon  x}{2} - i \cos \frac{ \varepsilon  x}{2}) &
32 i \sin \frac{\varepsilon x}{2}
\end{array}
\right).
$$
The first step of the asymptotic factorization leads to the
problem
\begin{equation}
\label{UnEx3_6} \Lambda^{+}(x) {{\hat N}}_{1,
\varepsilon }^{+}(x) + {{\hat N}}_{1, \varepsilon
}^{-}(x) \Lambda^{-}(x) = \hat{M}_{0, \varepsilon }(x),
\end{equation}
where the matrix function $\hat{M}_{0, \varepsilon }(x)=\hat{N}_{
\varepsilon }(x)$ can be represented in the following form:
\begin{equation}
\label{UnEx3_8} \hat{M}_{0, \varepsilon }(x) = \hat{M}_{0,
\varepsilon }^{+}(x) + \hat{M}_{0, \varepsilon }^{-}(x),
\end{equation}
{and}
\begin{equation} \label{UnEx3_9} \hat{M}_{0, \varepsilon
}^{+}(x) = \left(
\begin{array}{cc}
\frac{- 8 i(1 - e^{- \varepsilon })}{x + i} + \frac{8 x i(e^{i  \varepsilon  x} - e^{-  \varepsilon })}{x^2 + 1} & \frac{12 i(1 - e^{- \varepsilon })}{x + i} + \frac{- 16 x i(e^{i  \varepsilon  x} - e^{-  \varepsilon })}{x^2 + 1} \\
\frac{- 6 i(1 - e^{- \varepsilon })}{x + i} + \frac{4 x i(e^{i
\varepsilon  x} - e^{-  \varepsilon })}{x^2 + 1} & \frac{8 i(1 -
e^{- \varepsilon })}{x + i} + \frac{- 8 x i(e^{i  \varepsilon  x} -
e^{-  \varepsilon })}{x^2 + 1}
\end{array}
\right),
\end{equation}
\begin{equation}
\label{UnEx3_10} \hat{M}_{0, \varepsilon }^{-}(x) = \left(
\begin{array}{cc}
\frac{- 8 i(1 - e^{- \varepsilon })}{x + i} + \frac{8 x i(e^{- i  \varepsilon  x} - e^{-  \varepsilon })}{x^2 + 1} & \frac{12 i(1 - e^{- \varepsilon })}{x + i} + \frac{- 8 x i(e^{- i  \varepsilon  x} - e^{-  \varepsilon })}{x^2 + 1} \\
\frac{- 6 i(1 - e^{- \varepsilon })}{x + i} + \frac{8 x i(e^{- i
\varepsilon  x} - e^{-  \varepsilon })}{x^2 + 1} & \frac{8 i(1 -
e^{- \varepsilon })}{x + i} + \frac{- 8 x i(e^{- i  \varepsilon  x}
- e^{-  \varepsilon })}{x^2 + 1}
\end{array}
\right).
\end{equation}

Bounded solutions to (\ref{UnEx3_6}) have to satisfy the relation
\begin{equation}
\label{UnEx3_11} {{\hat N}}_{1,\varepsilon}^{+}(x) =
(\Lambda^{+}(x))^{-1} \left[\hat{M}_{0,\varepsilon}^{+}(x) -
\hat{C}_0\right],\; {{\hat N}}_{1, \varepsilon}^{-}(x) =
\left[\hat{M}_{0,\varepsilon}^{-}(x) + \hat{C}_0\right]
(\Lambda^{-}(x))^{-1},
\end{equation}
where $\hat{C}_0 = (\hat{c}_ij^{0})$ is a constant matrix.

In this case
$$
\hat{m}_{0,12}^{+}(i) = 6(1 - e^{- \varepsilon }) - 8 \varepsilon
e^{- \varepsilon },\quad \hat{m}_{0,12}^{-}(-i) = - 6(1 - e^{-
\varepsilon }) + 4 \varepsilon  e^{- \varepsilon },
$$
and thus, the solvability condition $\hat{m}_{0,12}^{+}(i) =
-\hat{m}_{0,12}^{-}(-i)$ is satisfied only for $ \varepsilon  =
0$. For all $ \varepsilon  \not= 0$, there is no approximate
solution (up to $\varepsilon^2$) of the functional equation
(similarly to (\ref{UnEx2_7})).

\begin{remark}
\label{no_guided} Note that, by construction,
$$
\hat{G}_{\varepsilon}(x) - G_{\varepsilon}(x) = \left(\begin{array}{cc} 0 & O(\varepsilon) \\
0 & 0 \end{array}\right),\; \varepsilon\to 0.
$$
Hence, $\hat{G}_{\varepsilon}(x)$ presents an example of the
regular perturbation of $\hat{G}_0(x)$, for which no regular
k-guided perturbation ($k > 1$) exists while construction of a
singular perturbation remains an open problem.
\end{remark}

\section{Numerical examples and discussion}
\label{numerics}

In this section, we analyse the quality of the approximation
provided by the 2-guided perturbation performed in Section\,
\ref{preserve} \ref{Ex_preserve}.

First, we consider the case when $c_{21}=0$, and thus the limiting
value of the remainder, $\Delta K_{1\varepsilon}$ does not vanish
at infinity. Specifically, we estimate the element on the main
diagonal in the following way:
\[
\Delta k_{jj}(\infty) =16\left(1-e^{-\varepsilon}+\varepsilon e^{-\varepsilon}\right)=
64\varepsilon^2-96\varepsilon^3+O(\varepsilon^4),\quad \varepsilon \to0,\quad j=1,2.
\]

\begin{figure}[h!]
  \begin{center}
    \hspace{-8mm}\includegraphics [scale=0.35]{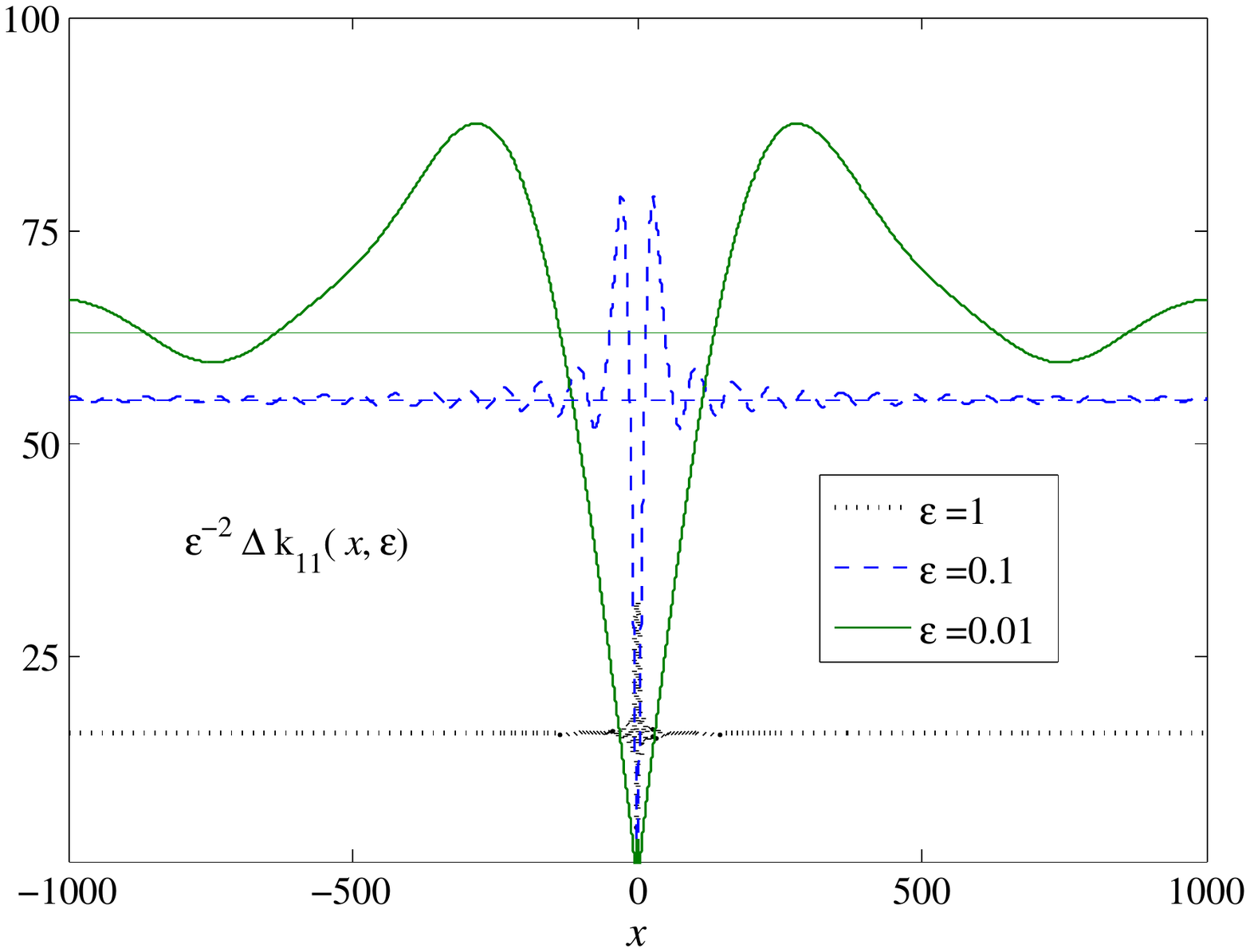}\hspace{-8mm}\includegraphics [scale=0.35]{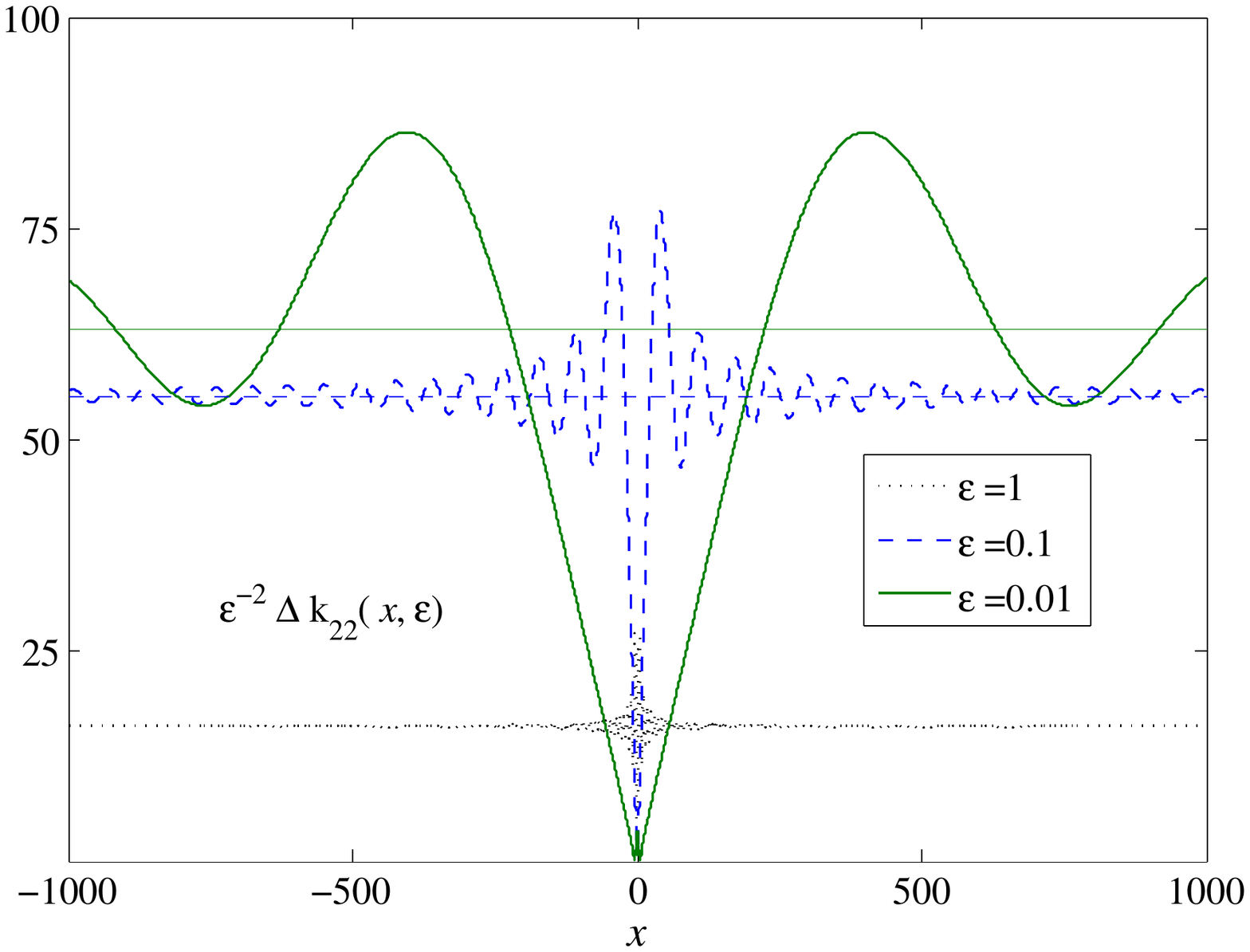}
       \end{center}
          \caption{\footnotesize{Diagonal elements, $\Delta k_{jj}(x,\varepsilon)$, $j=1,2$,
          of matrix $\Delta K_{1,\varepsilon}(x)$ defined in  for various values of parameter $\varepsilon$,
          the constant $c_{21}^{0}=0$.
    The elements are normalised to the value of parameter $\varepsilon^2$.}}
\label{f1}
\end{figure}

In Fig. \ref{f1} a) and b), those components are presented in
their normalised forms. We can see that the estimate is true (see
the discussion on the small parameter following formula
(\ref{difference})). Furthermore, the matrix converges to its
limiting values more quickly for larger values of the small
parameter, while the oscillations decay more slowly for smaller
values.

\begin{figure}[h!]
  \begin{center}
     \hspace{-8mm}\includegraphics [scale=0.35]{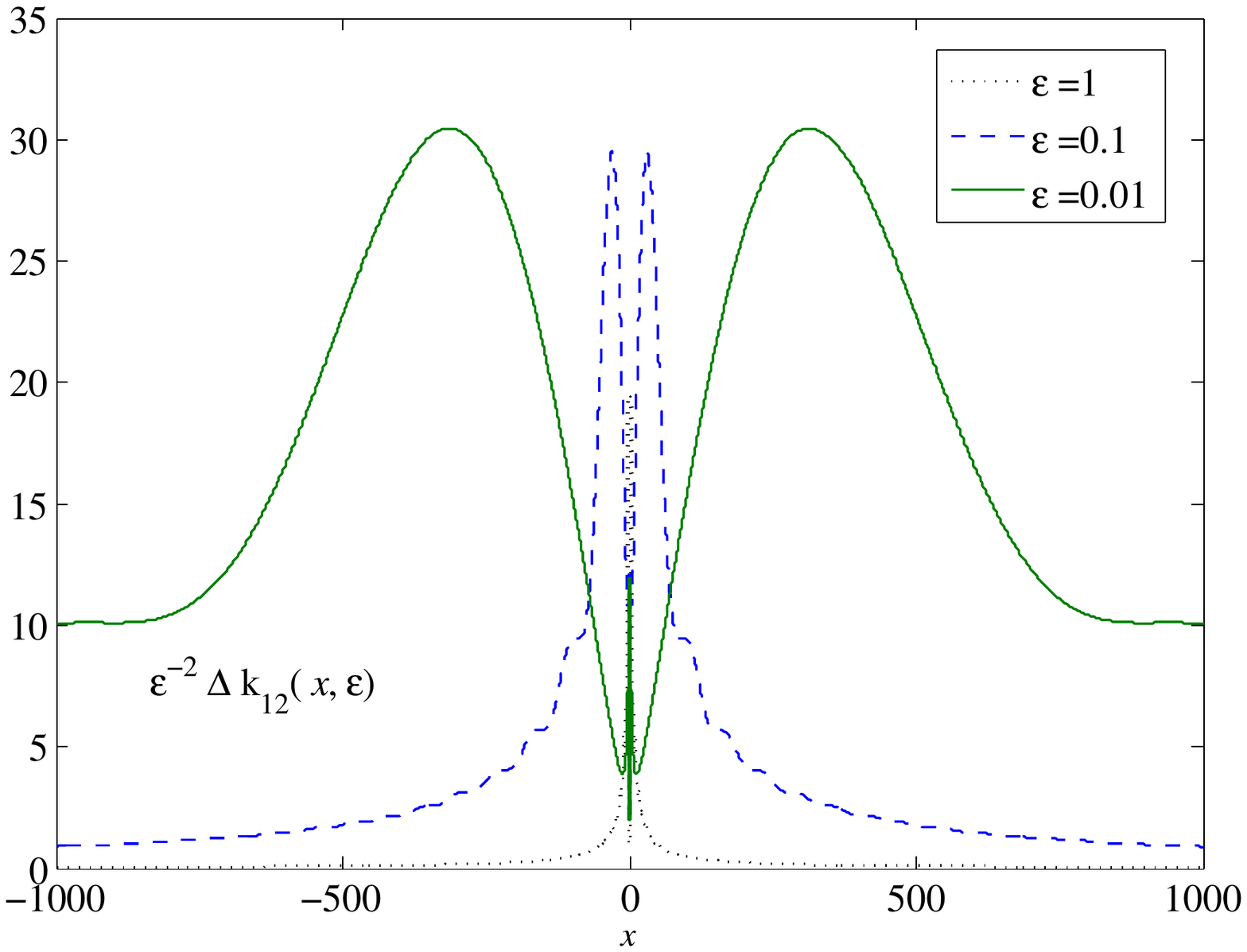}\hspace{-8mm}\includegraphics [scale=0.35]{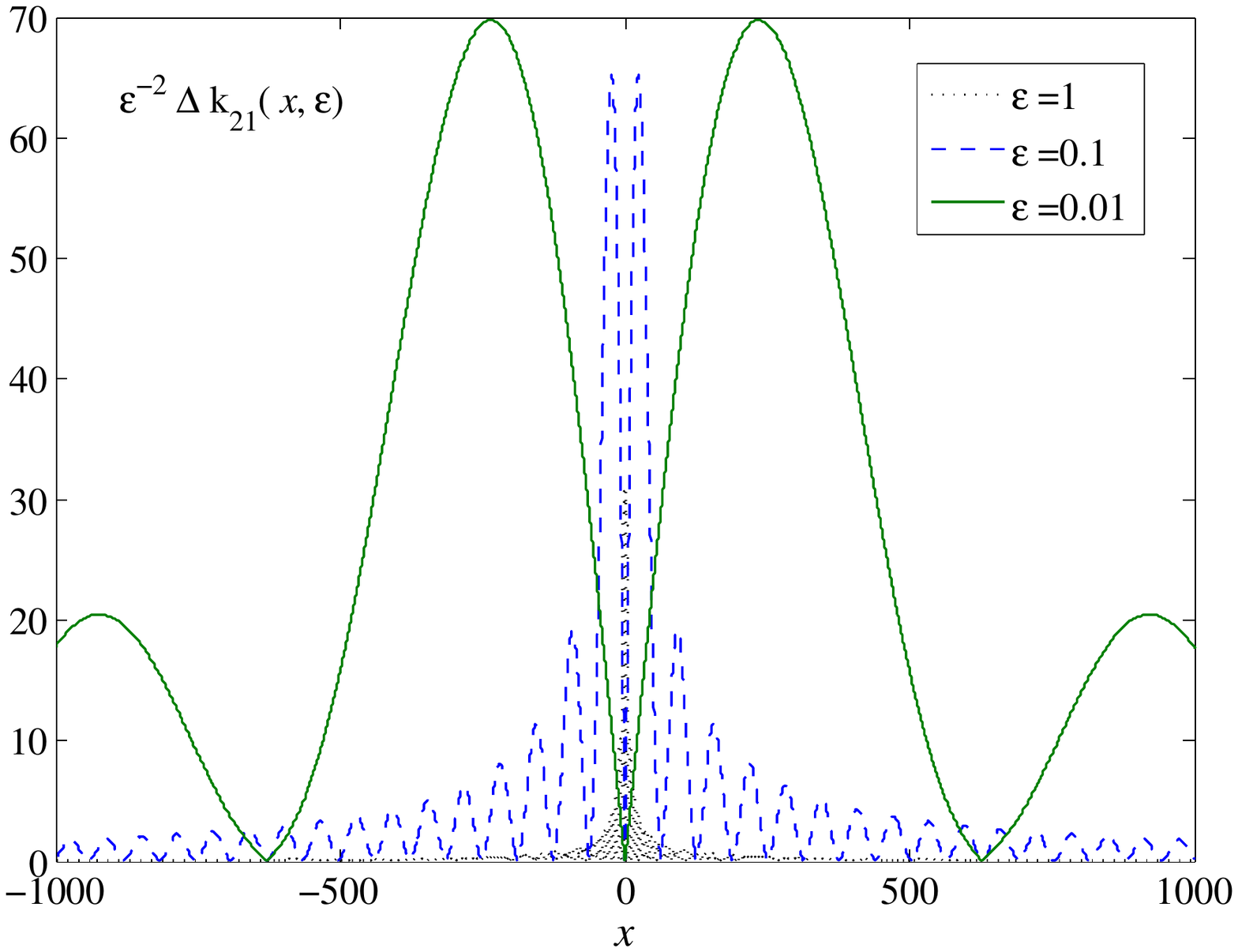}
         \end{center}
    \caption{\footnotesize{The other two elements, $\Delta k_{ij}(x,\varepsilon)$, $i+j=3$,
    of matrix $\Delta K_{1,\varepsilon}(x)$ {for} $\varepsilon=1;0.1;0.01$,
    {and} constant $c_{21}^{0}=0$. It is clear that both entries vanish at infinity $k_{ij}\to0$ as $|x|\to\infty$.}}
\label{f2}
\end{figure}

In Fig. \ref{f2} a) and b), the remaining two components are
depicted in the same normalised forms. Preserving the same
estimate, where $ \Delta  K_{1,\varepsilon}(x)=O(\varepsilon^2)$
as $\varepsilon\to0$, the components now decay to $O(x^{-1})$, as
$|x|\to\infty$. The trend is also clearly visible here, that the
smaller $\varepsilon$ is, the slower it converges to its limiting
value. In other words, the small parameter $\varepsilon$
determines the magnitude of the reminder matrix, but the
oscillations are larger in this case, and more pronounced along
the real axis.

Interestingly, the components on the main diagonal are comparable
in value, but not equal, while the remaining two differ in value
by almost a factor of two. Moreover, the latter are also two times
smaller in magnitude then the diagonal elements.

\begin{figure}[h!]
  \begin{center}
     \hspace{-8mm}\includegraphics [scale=0.35]{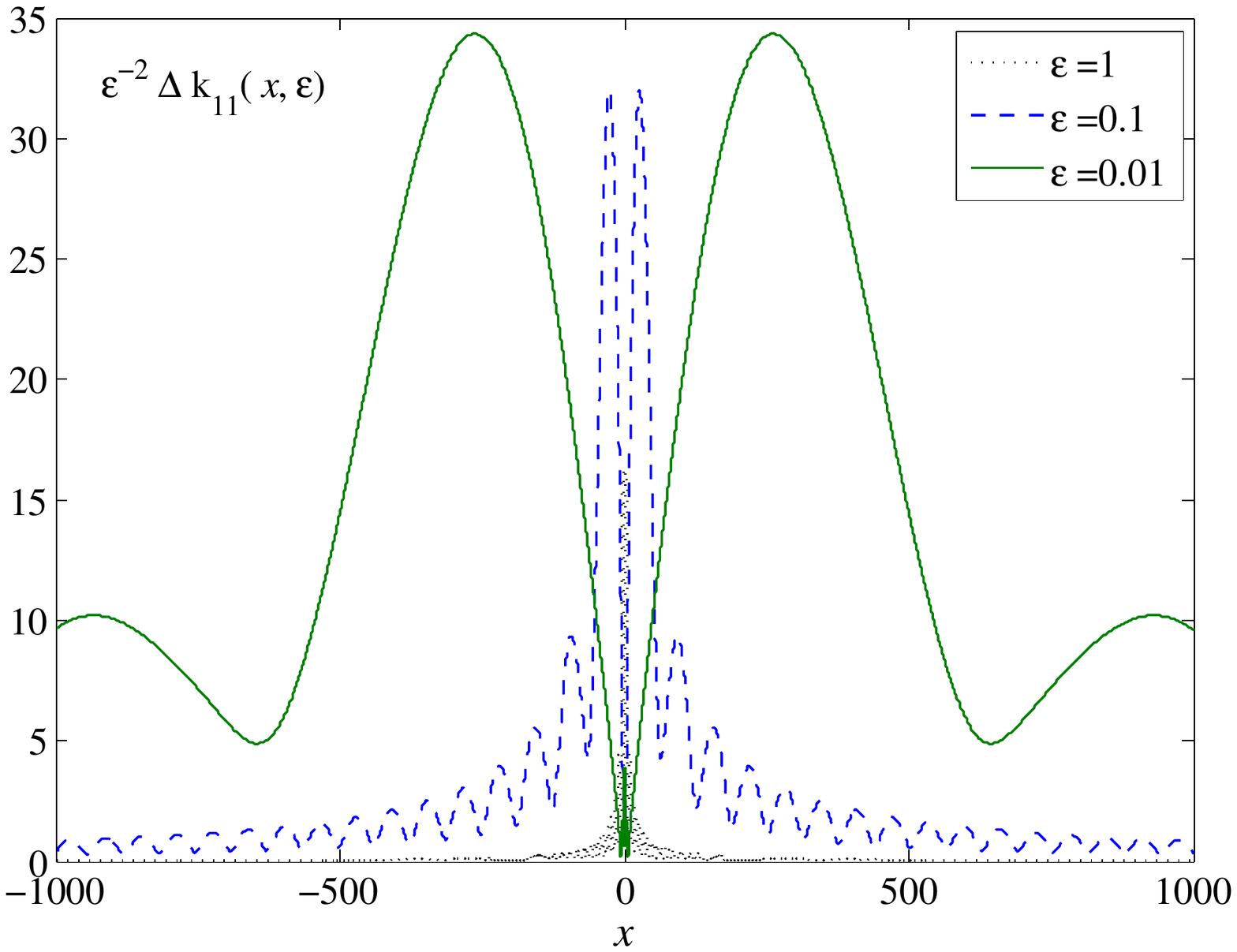}\hspace{-8mm}\includegraphics [scale=0.35]{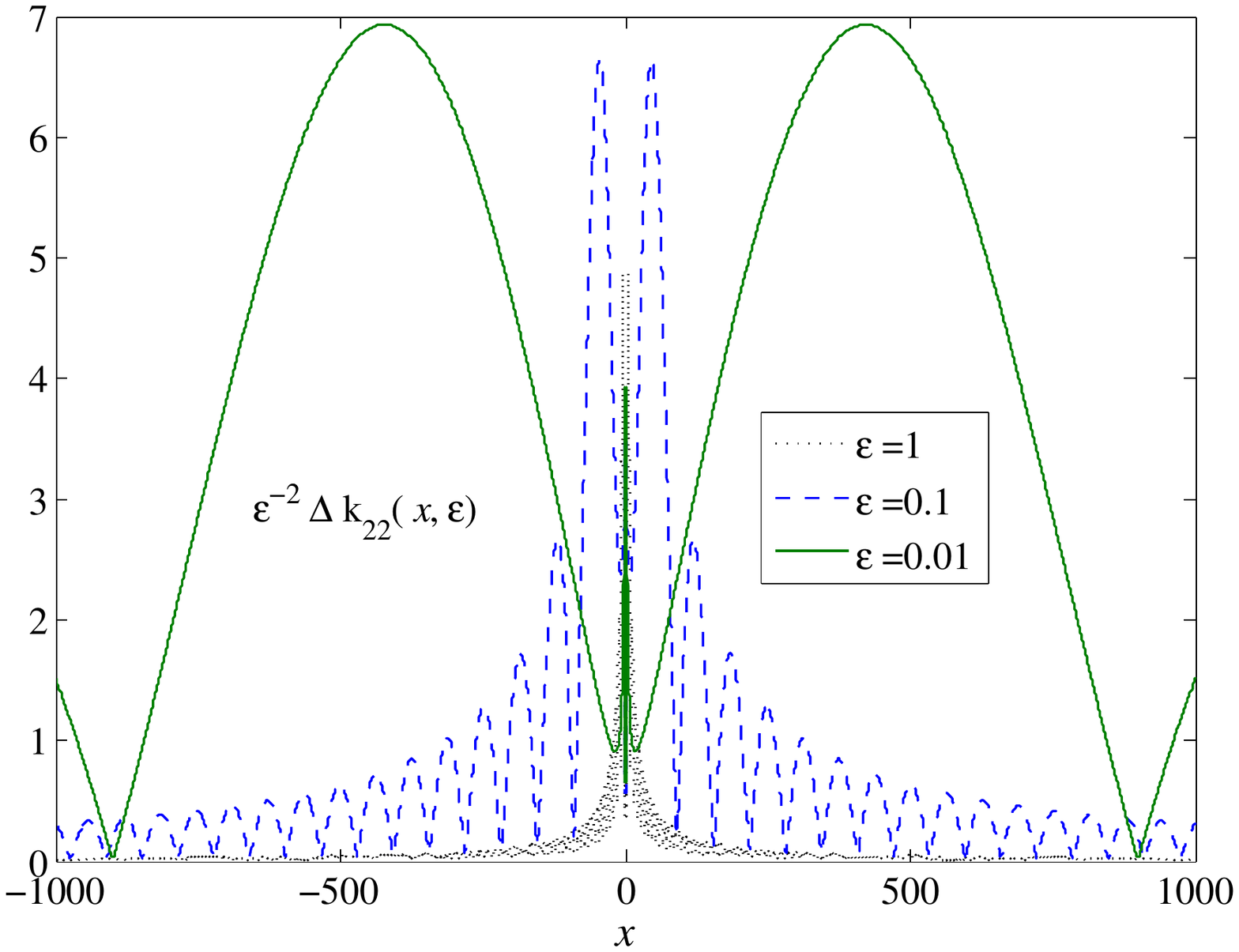}
       \end{center}
       \caption{\footnotesize{Diagonal elements, $\Delta k_{jj}(x,\varepsilon)$, $j=1,2$,
       of matrix $\Delta K_{1,\varepsilon}$, {for various} $\varepsilon$,
       {and} the constant $c_{21}^{0}=-8/3(1 - e^{-\varepsilon} + \varepsilon e^{-\varepsilon})$.
    The elements are normalised to the value of parameter $\varepsilon^2$. The horizontal lines show
    the limiting values of the normalised components at infinity.}}
    \label{f3}
\end{figure}

\begin{figure}[h!]
  \begin{center}
     \hspace{-8mm}\includegraphics [scale=0.35]{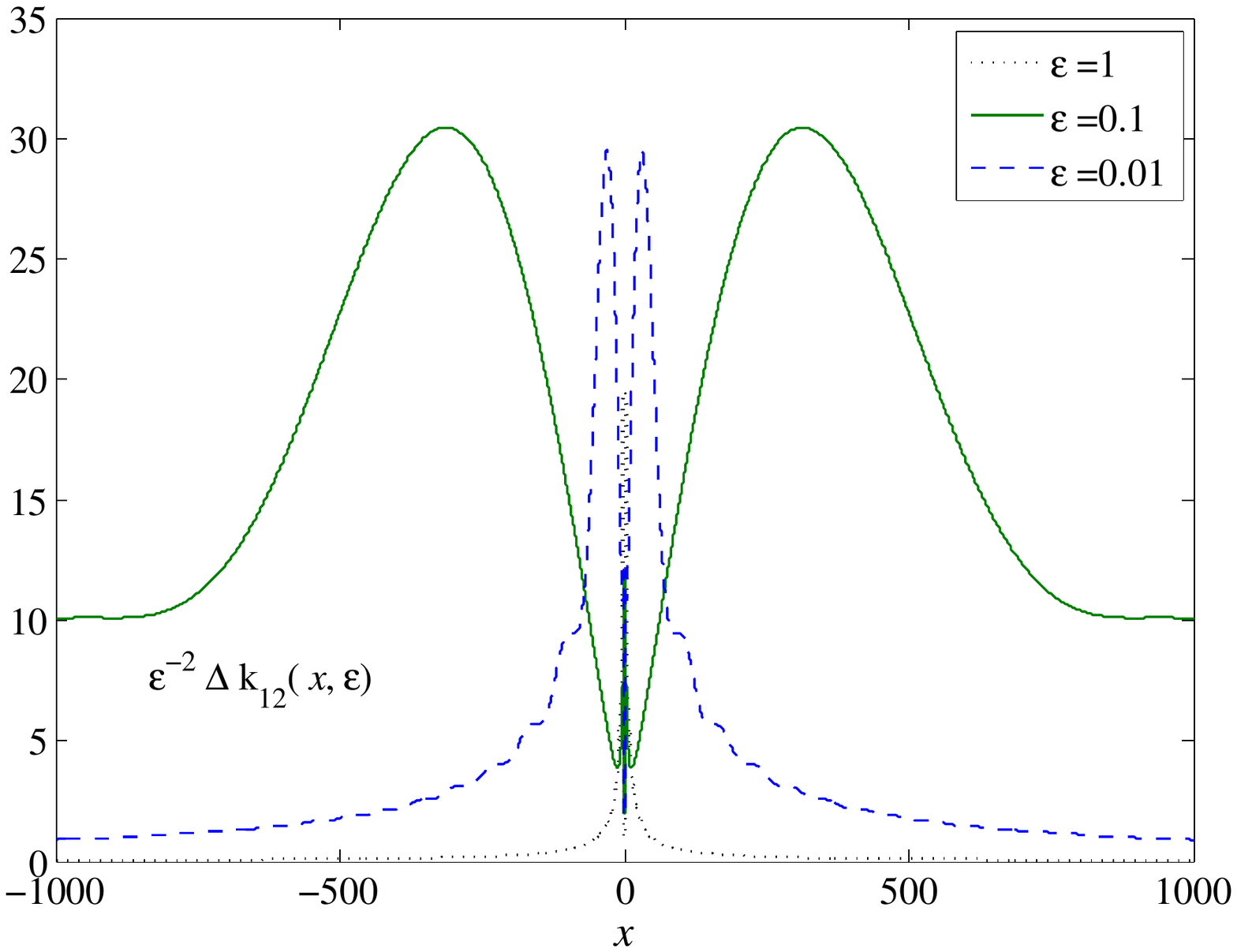}\hspace{-8mm} \includegraphics [scale=0.35]{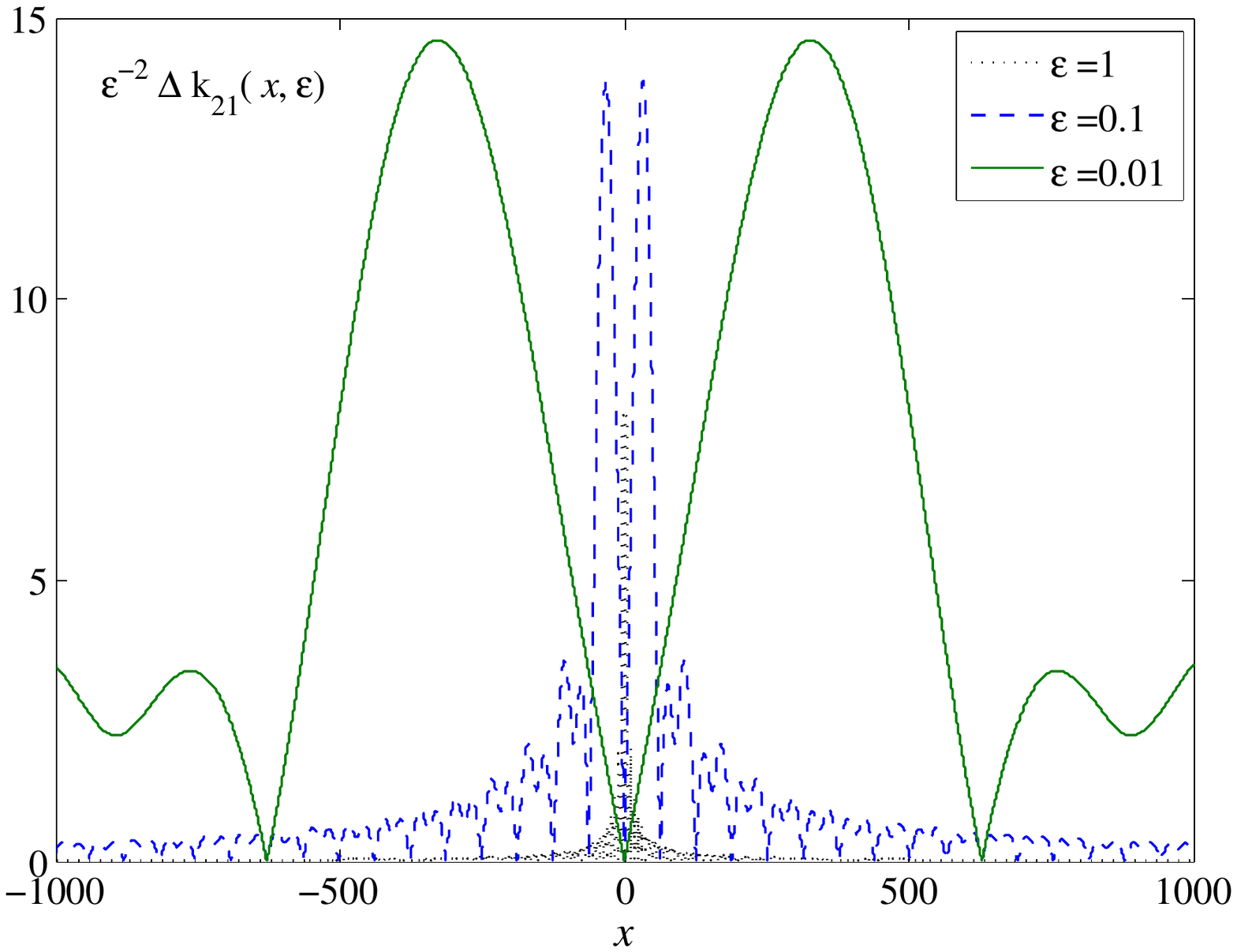}
         \end{center}
    \caption{\footnotesize{The other two elements, $\Delta k_{ij}(x,\varepsilon)$, $i+j=3$,
    of matrix $\Delta K_{1,\varepsilon}(x)$ {for} $\varepsilon=1;0.1;0.01$,
    {and} the constant $c_{21}^{0}=-8/3(1 - e^{-\varepsilon} + \varepsilon e^{-\varepsilon})$.
    It is clear that both entries vanish at infinity $k_{ij}\to0$ as $|x|\to\infty$.}}
    \label{f4}
\end{figure}

The situation changes when we consider the second case, where
$c_{21}^{0} =  - 8/3 (1 - e^{-\varepsilon} + \varepsilon
e^{-\varepsilon})$. The respective graphs are presented in Fig.
\ref{f3}, \ref{f4}. Now, all the components decay at infinity as
$O(x^{-1})$, as $|x|\to\infty$, and simultaneously have the same
estimate of ($O(\varepsilon^2)$) when $\varepsilon\to0$, as
predicted. The magnitudes of the components are, however, more
balanced in the sup norm $\|\Delta K_{1,\varepsilon}^{(1)}\|>2
\|\Delta K_{1,\varepsilon}^{(2)}\|$. This demonstrates that we can
choose an optimal approximation preserving some specified
requirement by varying the value of the arbitrary constant
$c_{21}$. Comparing these two cases, it is clear that the second
is preferable to the first for the reasons discussed above.

Any specific factorization will of course require its own
analyses. However, if the estimate
\[
G_\varepsilon(x)=o(1),\quad |x|\to\infty,
\]
is true, then the reminder can be estimated by
\[
\Delta K_{1,\varepsilon}(x)=
\left(
\begin{array}{cc}
 c_{11}^{0} &
c_{12}^{0} \\
c_{21}^{0} & c_{22}^{0}
\end{array}
\right)^2+o(1),\quad |x|\to\infty.
\]

We note that this property may change in the next step if we wish
to and can continue the approximation procedure, (the conditions
will remain valid for the next step). Here, the limiting values
for the first step will also play their role. We can deliver a
similar formula based on the two consequent approximations, where
two sets of constants will then be involved{:} $c_{jl}$ (first
step) and $d_{jl}$ (second step), $j,l=1,2$.

Judging by the magnitude of the reminder for both the presented
examples, we can conclude that the 2-guided perturbation may be
sufficient for practical purposes. Thus, if even one approximation
step is practically possible, meaning that the conditions
(\ref{sol_cond1})-(\ref{sol_cond5}) are satisfied, then we can use
this approximation directly in solving the Wiener-Hopf equation.

To close, we must highlight that, if conditions
(\ref{sol_cond1})-(\ref{sol_cond5}) for matrix $G_\varepsilon$,
with unstable partial indices, are not satisfied, the question of
how to compute a valuable approximate factorization for such a
matrix-function remains open.

{\bf Acknoledgements} The work is supported by the People Programme
(Marie Curie Actions) of the European Union's Seventh Framework Programme FP7/2007-
2013/ under REA grant agreement PIRSES-GA-2013-610547 - TAMER.
GM acknowledges support from a Royal Society Wolfson Research Merit Award.


\end{document}